\documentclass[11pt]{article}

\usepackage[english]{babel}
\usepackage[latin1]{inputenc}
\usepackage{fancyhdr}
\usepackage{amsmath}
\usepackage{amssymb}
\usepackage{amsfonts}
\usepackage{amsthm}
\usepackage{latexsym}
\usepackage{natbib}
\usepackage{graphicx,subfigure}
\usepackage{geometry}

\newtheorem{theo}{Theorem}[section]
\newtheorem{lem}{Lemma}[section]
\newtheorem{cor}{Corollary}[section]

\newtheorem{prop}{Proposition}[section]
\newtheorem{fact}{Fact}
\newtheorem{rem}{Remark}[section]

\begin{document}
\newcommand{\1}{{\rm 1}\kern-0.26em{\rm I}}
\newcommand{\E}{{\rm I}\kern-0.22em{\rm E}}
\renewcommand{\P}{{\rm I}\kern-0.22em{\rm P}}
\newcommand{\R}{{\rm I}\kern-0.20em{\rm R}}
\newcommand{\Anun}{A_{n,1}}
\newcommand{\alp}{\alpha}
\newcommand{\alpn}{\alpha_n}
\newcommand{\anstar}{a^{\star}_n}
\newcommand{\apcr}{\emph{apcr}}
\newcommand{\bnstar}{b^{\star}_n}
\newcommand{\BB}{{\cal B}}
\newcommand{\BBB}{\mathbb{B}}
\newcommand{\betn}{\beta_n}
\newcommand{\bx}{\boldsymbol{x}}
\newcommand{\bX}{\boldsymbol{X}}
\newcommand{\cadlag}{\emph{cadlag}}
\newcommand{\CC}{\mathcal{C}}
\newcommand{\CCC}{\mathbb{C}}
\newcommand{\DD}{{\mathcal{D}}}
\newcommand{\del}{\delta}
\newcommand{\deli}{\delta_i}
\newcommand{\delin}{\delta_{i,n}}
\newcommand{\demo}{\noindent {\textsc{Proof}.}}
\newcommand{\EE}{\mathcal{E}}
\newcommand{\EEE}{\mathbb{E}}
\newcommand{\EEnstar}{{\mathcal{E}}^{\star}_n}
\newcommand{\epsilone}{\varepsilon}
\newcommand{\e}{\varepsilon}
\newcommand{\epstendzer}{\varepsilon \downarrow 0}
\newcommand{\etanstar}{\eta^{\star}_n}
\newcommand{\etanstarun}{\eta^{\star (1)}_n}
\newcommand{\etanNstarun}{\eta^{\star (1)}_{n,N}}
\newcommand{\findem}{$\sqcup\!\!\!\!\sqcap$ \\}
\newcommand{\fcond}{\widehat{f}_{n,h,\ell}^{\star}}
\newcommand{\fcondn}{\widehat{f}_{n,h,\ell_n}^{\star}}
\newcommand{\FF}{{\mathcal{F}}}
\newcommand{\Fnstar}{F^{\star}_n}
\newcommand{\Fnhtstar}{\widetilde{F}^{\star}_{n,h}}
\newcommand{\Fnhhstar}{\widehat{F}^{\star}_{n,h}}
\newcommand{\Fnstarmoins}{F^{\star}_{n-}}
\newcommand{\fnstar}{f^{\star}_n}
\newcommand{\fnhstar}{f^{\star}_{n,h}}
\newcommand{\fnhnstar}{f^{\star}_{n,h_n}}
\newcommand{\gam}{\gamma}
\newcommand{\Gam}{\Gamma}
\newcommand{\Gamklj}{\Gamma_{k,l,j}}
\newcommand{\GG}{\mathcal{G}}
\newcommand{\Gnmoins}{G_{n-}}
\newcommand{\Gnstar}{G^{\star}_n}
\newcommand{\Gnstarm}{G^{\star}_{n-}}
\newcommand{\Gnstarbar}{\overline{G}^{\star}_n}
\newcommand{\Gnstarmoins}{G^{\star}_{n-}}
\newcommand{\HHnj}{{\cal{H}}^{(j)}_n}
\newcommand{\HHnun}{{\cal{H}}^{(1)}_n}
\newcommand{\HHnze}{{\cal{H}}^{(0)}_n}
\newcommand{\hd}{h^{1/d}}
\newcommand{\hk}{h'_{n_k}}
\newcommand{\hkl}{h'_{n_k, l}}
\newcommand{\hkld}{h^{\prime\ d}_{n_k, l}}
\newcommand{\hj}{H^{(j)}}
\newcommand{\hmoinsze}{H_-^{(0)}}
\newcommand{\hmoinsun}{H_-^{(1)}}
\newcommand{\hnj}{H_n^{(j)}}
\newcommand{\hnmoins}{H_{n-}}
\newcommand{\hnmoinsun}{H_{n-}^{(1)}}
\newcommand{\hnmoinsze}{H_{n-}^{(0)}}
\newcommand{\hnun}{H_n^{(1)}}
\newcommand{\hnze}{H_n^{(0)}}
\newcommand{\hpn}{h'_n}
\newcommand{\hsn}{h''_n}
\newcommand{\hun}{H^{(1)}}
\newcommand{\hze}{H^{(0)}}
\newcommand{\I}{I_d}
\newcommand{\infh}{\inf_{h\in [h'_n,h''_n]}}
\newcommand{\infk}{\inf_{k\in [k'_n,k''_n]}}
\newcommand{\inth}{h\in [h'_n,h''_n]}
\newcommand{\inthk}{h\in [h'_{n_k},h''_{n_k-1}]}
\newcommand{\inthkl}{h\in [h'_{n_k,l},h'_{n_k,l+1}]}
\newcommand{\intk}{k\in [k'_n,k''_n]}
\newcommand{\intjjl}{1\leq j\leq J_l}
\newcommand{\intlrk}{0\leq l\leq R_k}
\newcommand{\intlrkm}{0\leq l\leq R_k-1}
\newcommand{\intox}{\int_0^x}
\newcommand{{\iunn}}{1\leq i \leq n}
\newcommand{\KK}{{\mathcal{K}}}
\newcommand{\KKnstar}{{\mathcal{K}}^{\star}_n}
\newcommand{\KKK}{\mathbb{K}}
\newcommand{\KM}{Kaplan-Meier}
\newcommand{\lab}{\lambda}
\newcommand{\lac}{\bigg\{}
\newcommand{\Lac}{\Big\{}
\newcommand{\lambdanhstar}{\lambda^{\star}_{n,h}}
\newcommand{\lcro}{\bigg[}
\newcommand{\Lcro}{\Big[}
\newcommand{\LCro}{\Bigg[}
\newcommand{\ld}{l=1, \cdots, d}
\newcommand{\limk}{\lim_{k\rightarrow\infty}\;}
\newcommand{\limn}{\lim_{n\rightarrow\infty}\;}
\newcommand{\LL}{{\cal L}}
\newcommand{\LLnstar}{{\mathcal{L}}^{\star}_n}
\newcommand{\lpar}{\bigg(}
\newcommand{\Lpar}{\Big(}
\newcommand{\maxn}{\max_{n\in N_k}}
\newcommand{\maxjjl}{\max_{1\leq j \leq J_l}}
\newcommand{\maxlrk}{\max_{0\leq l \leq R_k}}
\newcommand{\maxlrkm}{\max_{0\leq l \leq R_k-1}}
\newcommand{\maxoscillolj}{\max_{\substack{0\leq l \leq R_k-1\\1\leq j \leq J_l}}}
\newcommand{\MM}{{\mathcal{M}}}
\newcommand{\mnhpsistarh}{\widehat{m}^{\star}_{\psi,n,h}}
\newcommand{\n}{n^{1/2}}
\newcommand{\ninf}{n\rightarrow\infty}
\newcommand{\NNN}{\mathbb{N}}
\newcommand{\NN}{\mathcal{N}}
\newcommand{\nono}{\nonumber}
\newcommand{\OO}{{\mathcal{O}}}
\newcommand{\p}{{\rm I}\kern-0.18em{\rm P}}
\newcommand{\PPP}{ \mathbb{P}}
\renewcommand{\proof}{\noindent {\bf \textsc{Proof}.}}
\newcommand{\Q}{\mathbb{Q}}
\newcommand{\qun}{Q^{(1)}}
\newcommand{\qze}{Q^{(0)}}
\newcommand{\rac}{\bigg\}}
\newcommand{\Rac}{\Big\}}
\newcommand{\rcro}{\bigg]}
\newcommand{\Rcro}{\Big]}
\newcommand{\RCro}{\Bigg]}
\newcommand{\rpar}{\bigg)}
\newcommand{\Rpar}{\Big)}
\newcommand{\RRR}{\mathbb{R}}
\newcommand{\sigmapsi}{\sigma(\psi)}
\newcommand{\SSS}{\mathbb{S}}
\newcommand{\sumin}{\sum_{i=1}^n}
\newcommand{\supab}{\sup_{a\leq t \leq b}}
\newcommand{\supabx}{\sup_{a\leq x \leq b}}
\newcommand{\supabpx}{\sup_{a\leq x \leq b'}}
\newcommand{\supabp}{\sup_{a\leq t \leq b'}}
\newcommand{\suphzeun}{\sup_{\substack{0\leq s,t \leq 1\\ |t-s|\leq
h}}}
\newcommand{\suphabp}{\sup_{\substack{a\leq s,t \leq b'\\ |t-s|\leq
h}}}
\newcommand{\suph}{\sup_{h\in [h'_n,h''_n]}}
\newcommand{\suphk}{\sup_{h\in [h'_{n_k},h''_{n_{k-1}}]}}
\newcommand{\suphkl}{\sup_{h\in [h'_{n_k,l},h'_{n_k,l+1}]}}
\newcommand{\suphx}{\sup_{\substack{h\in [h'_n,h''_n] \\ \bx\in I}}}
\newcommand{\supk}{\sup_{k\in [k'_n,k''_n]}}
\newcommand{\suposcillohx}{\sup_{\substack{\boldsymbol{x}\in \Gamma_{k,l,j} \\h\in [h'_{n_k,l},h'_{n_k,l+1}]}}}
\newcommand{\suppsi}{\sup_{\psi\in\FF}}
\newcommand{\supxGam}{\sup_{x\in \Gamma_{k,l,j}}}
\newcommand{\sups}{\sup_{0\leq s \leq 1}}
\newcommand{\supx}{\sup_{x\in I}}
\newcommand{\tend}{\rightarrow}
\newcommand{\toh}{\tau_h}
\newcommand{\tohm}{\tau^{-1}_h}
\newcommand{\tohpm}{\tau^{-1}_{h'_n}}
\newcommand{\UU}{{\cal U}}
\newcommand{\unifh}{h\in [h'_n,h''_n]}
\newcommand{\vpinun}{\varpi_n^{(1)}}
\newcommand{\wh}{\widehat}
\newcommand{\wnun}{w_n^{(1)}}
\newcommand{\xin}{\xi_n}
\newcommand{\xklj}{{\boldsymbol x}_{k,l,j}}
\newcommand{\XX}{{\mathcal{X}}}
\newcommand{\zetan}{\zeta_n}
\newcommand{\Zidelin}{\{(Z_i, \delta_i):1\leq i \leq n\}}
\newcommand{\Zin}{\Z_{i,n}}
%\begin{frontmatter}

\title{{\bf Uniform limit laws of the logarithm for nonparametric estimators of the regression function in presence of censored data.}}
\author{MAILLOT, B. $^{(1)}$ and VIALLON, V.$^{(1,2),}$\thanks{Correspondence to : Vivian Viallon, Laboratoire de Biostatistique, Hôpital Cochin, Faculté de
Médecine, Université Paris-Descartes, 24 rue du Faubourg Saint
Jacques, 75014 Paris, France. mail :
vivian.viallon@univ-paris5.fr}}

\date{{\small $^{(1)}$ Laboratoire de Statistique Théorique et Appliquée (LSTA), Université Paris 6, 175, rue du Chevaleret, 75013
Paris,
France\\
$^{(2)}$ Laboratoire de Biostatistique, Hôpital Cochin, Faculté de
Médecine, Université Paris-Descartes, 24 rue du Faubourg Saint
Jacques, 75014 Paris, France. }}
\maketitle
\begin{abstract}
In this paper, we establish uniform-in-bandwidth limit laws of the
logarithm for nonparametric Inverse Probability of Censoring
Weighted (I.P.C.W.) estimators of the multivariate regression
function under random censorship. A similar result is deduced for
estimators of the conditional distribution function. The
uniform-in-bandwidth consistency for estimators of the conditional
density and the conditional hazard rate functions are also derived
from our main result. Moreover, the logarithm laws we establish
are shown to yield almost sure simultaneous asymptotic confidence
bands for the functions we consider. Examples of confidence bands
obtained from simulated data are displayed.\vskip5pt

\vskip5pt

\end{abstract}

%\begin{keyword}
%\kwd{Estimation non paramétrique} \kwd{Données censurées}
%\kwd{Modèle additif} \kwd{Intégration marginale} \kwd{Régression
%multivariée.}
%\end{keyword}
%
%\begin{keyword}[class=AMS]
%\kwd{62N01} \kwd{62G08.}
%\end{keyword}
%
%\end{frontmatter}

\section{Introduction-Motivations}
\setcounter{equation}{0} \setcounter{theo}{0} \setcounter{lem}{0}
\setcounter{cor}{0} \setcounter{rem}{0} \setcounter{prop}{0}
\setcounter{fact}{0}

\noindent Nonparametric estimators of functionals of the
conditional law (such as the regression function or the
conditional distribution function) are known to provide a suitable
and efficient means to catch the possibly complex relation between
a given variable of interest and some explanatory covariates.
Because of this obvious practical interest, many authors have
studied the (asymptotic) properties of such estimators (see, e.g.,
\cite{BOSLEC87}, \cite{Hardle90}, \cite{Gyorfietal02}). Fewer
works deal with the special case where the variable of interest is
censored (\cite{FanGijbels}, \cite{Stute99}). Yet, this situation
arises in many statistical applications, as medical research,
reliability, ... and it is therefore of paramount importance to
build and study estimators adapted to the censored setting. When
the variable of interest is subject to right-censoring,
transformations of the observed data are generally needed to
derive inference on the underlying (conditional) distribution (see
\cite{BuckleyJames}). Estimates based on these transformations are
usually referred to as \emph{synthetic data} estimates in the
literature. In the case of the regression function estimation,
Fan and Gijbels \cite{FanGijbels} especially proposed a transformation leading to
a local version of the Buckley-James estimator. In this paper, we
make use of an alternative transformation which leads to
\emph{Inverse Probability of Censoring Weighted} $[I.P.C.W.]$
estimators. I.P.C.W. type estimators have recently gained
popularity in the censored literature. To our mind, they basically
present two particularly appealing properties. First, and as it
will be especially made clear in the proofs of our forthcoming
results (see also \cite{Carbonez1995}, \cite{Kohler2002},
\cite{BrunelComte1}, \cite{Kohler2003} and the references
therein), their asymptotic behavior can be easily deduced from
that of analogous estimators in the uncensored case. Second, their
computation is straightforward. In that sense, they are appealing
for both theoretic and applied statistics purposes. It is however
noteworthy that methodology we propose here for $I.P.C.W$-type
estimates shall apply with minor modifications to cope with other
synthetic data estimates.\vskip5pt

\noindent The present paper is organized as follows. First, we
introduce the main notations and hypotheses needed for our task.
Then, following the methodology developed  in the uncensored case
by Einmahl and Mason \cite{EM00}, we establish a uniform-in-bandwidth law of the
logarithm for a nonparametric I.P.C.W. estimator of the regression
function (see Theorem \ref{theo_carbo} below). This result
corresponds to the almost sure and uniform-in-bandwidth version of
Theorem 3.1 in \cite{Viallon_CRAS1}. At this point we shall stress
the reader attention on the fact that, as was first shown by Deheuvels and Mason
\cite{DM04} (see also \cite{EM05}, \cite{DonyEM} and the relevant
references therein), such uniform-in-bandwidth limit laws turn out
to be of particular interest in practice because they allow for
establishing uniform consistency of data-driven (and then random)
bandwidth estimators. In Section \ref{section_corollaires} we
derive a similar law of the logarithm for an estimator of the
conditional distribution function and we establish the
uniform-in-bandwidth consistency for estimators of the conditional
density and the conditional hazard rate functions. As was
especially pointed out by Deheuvels and Mason \cite{DM04} in the uncensored case,
limit laws of the logarithm provide themselves useful in the
construction of simultaneous confidence bands for the true
considered function. Such confidence bands, obtained from
simulated data, are displayed in Section \ref{section_IC}.
Finally, Section \ref{section_proofs} is devoted to the proofs of
our results.

\section{Notations and hypotheses}

\setcounter{equation}{0} \setcounter{theo}{0} \setcounter{lem}{0}
\setcounter{cor}{0} \setcounter{rem}{0} \setcounter{prop}{0}
\setcounter{fact}{0}

Consider a triple $({Y}, C, {\bf X})$ of random variables defined
in ${\R}\times {\R} \times {\R}^d$, $d\geq1$. Here $Y$ is the
variable of interest, $C$ a censoring variable and ${\bf
X}=(X_1,...,X_d)$ a vector of concomitant variables.  Throughout,
we work with a sample $\{(Y_i,C_i,{\bf X}_i)_{1\leq i \leq n}\}$
of independent and identically distributed replica of $(Y,C,{\bf
X})$, $n\geq1$. Actually, in the right censorship model, the pairs
$(Y_i,C_i)$, $1\leq i\leq n$ are not directly observed and the
corresponding information is given by $Z_{i}:=\min\{Y_{i},C_{i}\}$
and $\delta_{i}:=\1_{\{{Y_{i}} \leq C_i\}}$, $1\leq i\leq n$, with
$\1_E$ standing for the indicator function of $E$. Accordingly,
the observed sample is ${\cal D}_n =\{(Z_i, \deli,
\boldsymbol{X}_i), i=1, \ldots, n\}$. \vskip5pt

\noindent In the sequel, we impose the following assumptions upon
the distribution of $(\bX, Y)$. Denote by $I$ a given compact of
$\R^d$ with non empty interior and set, for any $\gamma>0$,
\[I^\gamma = \{\bx : \inf_{{\bf u}\in I}|\bx-{\bf u}|_{\R^d}\leq
\gamma\},\] with $|\cdot|_{\R^d}$ standing for the usual euclidian
norm on $\R^d$. We will assume that, for a given $\alpha>0$,
$(\bX, Y)$ [resp. $\bX$] has a density function $f_{\bX, Y}$
[resp. $f_{\bX}$] with respect to the Lebesgue measure on
$I^\alpha\times\R$ [resp. $I^\alpha$]. We will also assume that
the assumptions $(F.1$-$2)$ below hold. For $-\infty< t<\infty$,
set $F(t)=\P(Y\leq t)$, $G(t)=\P(C\leq t)$ and $H(t)=\P(Z\leq t)$,
the right continuous distribution functions of $Y$, $C$ and $Z$
respectively. For any right continuous distribution function $L$
defined on $\R$, further denote by $T_L=\sup\{t\in\R:L(t)<1\}$ the
upper point of the corresponding distribution.
\begin{tabbing}
$(F.1)\;\;$ \=\ For all ${\bx}\in I^{\alp}$,
$\displaystyle\lim_{{\bx}'\rightarrow {\bx};{\bx}'\in
I^{\alp}}f_{\boldsymbol{X},{Y}}({\bx}',{
y})=f_{\boldsymbol{X},{Y}}({\bx},{y})$ for almost every
${y}\leq T_H$. \\[0.2cm]
$(F.2)$ \>\ $f_{\boldsymbol{X}}$ is continuous and strictly
positive on $I^{\alp}$.
\end{tabbing}
Now consider a \emph{pointwise measurable} class  $\FF$ (see p.110
in \cite{vdVW96}) of real measurable functions defined on $\R$.
Throughout, $\FF$ will be assumed to form a \emph{VC subgraph
class} (see \S 2.6.2 in \cite{vdVW96}).\\
\noindent In this paper, we will mostly focus on the regression
function of $\psi(Y)$ evaluated at ${\bf X}={\bf x}$, for
$\psi\in\FF$ and $\bx\in I^\alpha$, given by
\begin{eqnarray}
m_{\psi}({\bf x}) =  \E\big(\psi (Y) \mid {\bf X}={\bf x}\big).
\label{fdereg}
\end{eqnarray}
To estimate $m_\psi$ when $Y$ is right-censored, the key idea of
I.P.C.W. estimators is as follows. Introduce the real valued
function $\Phi_\psi$ defined on $\R^2$ by
\begin{equation}
\Phi_\psi(y,c)= \frac{\1_{\{y\leq c\}}\psi(y\wedge c)}{1-
G(y\wedge c)}.\label{def_phipsi}
\end{equation}
Assuming the function $G$ to be known, first note that
$\Phi_\psi(Y_i,C_i)=\delta_i\psi(Z_i)/ (1-G(Z_i))$ is observed for
every $1\leq i\leq n$. Moreover, under the assumption $({\cal I})$
below,
\begin{tabbing}
$ ({\cal I})\;\;$ \=\ $C$ and $(Y,\boldsymbol{X})$ are
independent;
\end{tabbing}
we have
\begin{eqnarray}
m_{\Phi_\psi}({\bf x})& := &\E(\Phi_\psi(Y,C)|{\bf X}=\bx) \nonumber\\
& =& \E\Big\{\frac{\1_{\{Y\leq C\}}\psi(Z)}{1-{G}(Z)}\big|\ {\bf X}=\bx\Big\}\nonumber\\
& = &\E \Big\{\frac{\psi(Y)}{1-{G}(Y)} \E\big(\1_{\{Y\leq
C\}}|{\bf X}, Y\big)\big|\ {\bf X}=\bx\Big\}\nonumber\\
& = &m_{\psi}({\bf x}).\label{key_idea}
\end{eqnarray}
Therefore, any estimate for $m_{\Phi_\psi}$, which can be built on
fully observed data, turns out to be an estimate for $m_\psi$ too.
Thanks to this property, most statistical procedures known to
provide estimates of the regression function in the uncensored
case can be naturally extended to the censored case. For instance,
kernel-type
estimates are particularly easy to construct. %\begin{rem}
%Assumption $({\cal I})$ is stronger than the conditional
%independence between $Y$ and $C$ given $\bX$, which is generally
%made in the censored regression field. However, $({\cal I})$ is
%essentially needed when I.P.C.W type estimates are involved, as
%can be seen through (\ref{key_idea}) above. The case where $({\cal
%I})$ does not hold will be addressed in Section
%\ref{Res_Indep_cond}.
%\end{rem}
Let $K$ be a kernel function defined on $\R^d$, that is a
measurable function such that $\int_{\R^d}K(\bx)d\bx=1$, and set,
for ${\bf x}\in I$, $h>0$, $1\leq i\leq n$,
\begin{equation}\label{poids_carbo}
\varpi_{n,h,i}({\bx}) := K\Big(
\frac{{\bx}-\boldsymbol{X}_i}{h}\Big) \Big/ \sum_{j=1}^nK\Big(
\frac{{\bx}-\boldsymbol{X}_j}{h}\Big).
\end{equation}
\noindent Unless otherwise specified, we will let $h>0$ vary in
such a way that $h'_n\leq h\leq h''_n$, where $\{h'_n\}_{n\geq1}$
and $\{h''_n\}_{n\geq1}$ are two sequences of positive constants
such that $0<h'_n\leq h''_n<\infty$ and, for either choice of
$h_n=h'_n$ or $h_n=h''_n$, the conditions $(H.1$-$2$-$3)$ below
are fulfilled by $\{h_n\}_{n\geq1}$.
\begin{tabbing}
$(H.1)\;\;\;\;$ \=  $h_n \downarrow 0$, $0<h_n<1$, and $nh_n^d\uparrow\infty$;\\[0.2cm]
$(H.2)$ \> $nh_n^d/\log n \rightarrow \infty$ as $n
\rightarrow \infty$;\\[0.2cm]
$(H.3)$ \> $\log(1/h_n)/\log\log n\rightarrow \infty$ as $n
\rightarrow \infty$.
\end{tabbing}
At this point, we shall make the reader note that, under
$(H.1$-$2$-$3)$ and $(F.2)$, the denominator involved in the
expression of the functions $\varpi_{n,h,i}$ is almost surely
strictly positive on $I$ for $n$ large enough, for every $1\leq
i\leq n$ and all $\inth$.\vskip5pt

\noindent In view of (\ref{def_phipsi}), (\ref{key_idea}) and
(\ref{poids_carbo}), whenever $G$ is known, a kernel estimator of
$m_{\psi}({\bf x})$ is given by
\begin{equation}
\widehat{m}_{\psi,n,h}({\bx}) := \sum_{i=1}^n
\varpi_{n,h,i}({\bx})\frac{\deli \psi(Z_i)}{1-G(Z_i)}.
\label{estcarboGconnu}
\end{equation}
\noindent In practice however, the function $G$ is generally
unknown and then has to be estimated. We will denote by $\Gnstar$
the Kaplan-Meier estimator of the function
$G$ \cite{KaplanMeier}. Namely, adopting the conventions $\prod_{\varnothing}=1$ and
$0^0=1$ and setting $N_n(x)=\sum_{i=1}^n\1_{\{Z_i\geq x\}}$, we
have
\begin{eqnarray}
\Gnstar(u) = 1-\prod_{i:Z_i\leq
u}{\Big(\frac{N_n(Z_i)-1}{N_n(Z_i)}\Big)}^{(1-\delta_i)},\mbox{
for all }u\in{\R}. \label{Gnstar}
\end{eqnarray}
Given these notations,  the following estimator of $m_\psi(\bx)$
can be proposed (see \cite{Kohler2002} for instance),
\begin{equation}
\widehat{m}_{\psi,n,h}^{\star}({\bx}) := \sum_{i=1}^n
\varpi_{n,h,i}({\bx})\frac{\deli \psi(Z_i)}{1-\Gnstar(Z_i)}.
\label{estcarbo1}
\end{equation}
Adopting the convention $0/0=0$, this quantity is properly defined
since $\Gnstar(Z_i)=1$ if and only if $Z_i=Z_{(n)}$ and
$\delta_{(n)}=0$, where $Z_{(k)}$ is the $k$-th ordered statistic
associated to the sample $(Z_1,...,Z_n)$ for $k=1,...,n$ and
$\delta_{(k)}$ is the $\delta_j$ corresponding to
$Z_{(k)}=Z_j$.\vskip5pt

\noindent As mentioned in \cite{GrossLai}, functionals of the
(conditional) law can generally not be estimated on the complete
support when the variable of interest is right-censored.
Accordingly, in order to establish our results, we will work under
the assumption $({\bf A})$ that will be said to hold if either
$({\bf A})(i)$ or $({\bf A})(ii)$ below holds.
\begin{tabbing}
$({\bf A})(i)\;\;$ \=\ There exists
a $\tau_0<T_H$ such that, for all $\psi\in\FF$,  $\psi=0$ on $(\tau_0,\infty)$.\\[0.2cm]
$({\bf A})(ii)$ \>\ $(a)\;\;$ For a given $0<p\leq 1/2$,
$\int_0^{T_H}
(1-F)^{-p/(1-p)}dG<\infty$;\\
\>\ $(b)\;\;$ $T_F<T_G$ and $(Y,C)\in\R^+\times\R^+;$\\
\>\ $(c)\;\; n^{2p-1}h^{''-d}_{n}|\log h''_{n}|\rightarrow
\infty$, as $n\rightarrow\infty$.
\end{tabbing}
\begin{rem} The assumption $({\bf A})(ii)$ will be needed in
our proofs when considering the estimation of the "classical"
regression function, which corresponds to the choice $\psi(y)=y$.
On the other hand, rates of convergence for estimators of
functionals such as the conditional distribution function
$\p(Y\leq t|{\bf X})$ can be obtained under weaker conditions,
when restricting ourselves to $t\in[-\infty,\tau_0]$ with
$\tau_0<T_H$. Compare the conditions in Theorems \ref{theo_carbo}
and Corollary \ref{cor_fdr} below.
\end{rem}

\noindent Besides the above assumptions, we will impose the
following additional hypotheses to establish our main results.

\begin{tabbing}
$(F.3)\;\; $ \=\  The variable $C$ has a Lebesgue density function
$f_C$ on $(-\infty, T_H]$.\\
$(F.4)$ \>\ The class of functions $\FF$ is bounded (in the sense
that $\FF$ has a measurable and \\ \>\ uniformly bounded envelope
function $\Upsilon({y})\geq\sup_{\psi\in\FF}\psi({
y})$, ${y}\leq T_H$).\\
$(F.5)$ \>\ The class of functions $\MM:=\{m_{\psi
}/f_{\boldsymbol{X}}, \psi\in\FF\}$ is relatively compact with
 \\ \>\ respect to the sup-norm  topology on $I^\alpha$.
\end{tabbing}

\noindent It is noteworthy that $(F.4$-$5)$ are automatically
fulfilled for the particular choice $\FF=\{\1_{(-\infty,t]},t\leq
T_H\}$ (see pp. 6-7 in \cite{EM00}). This property will enable us
to easily describe uniform consistency for estimators of the
conditional distribution function $F(t;\cdot):=\P(Y\leq t
|\bX=\cdot)$ over $t\in(-\infty,T_H)$(see Corollary \ref{cor_fdr}
in Section \ref{section_corollaires}). \vskip5pt

\noindent Turning our attention to the kernel $K$, set
$\KK:=\{K(\lambda(\cdot-{\bf y})), {\bf y}\in\R^d, \lambda>0\}$
and denote by $\mathcal{N}(\e,\mathcal{K})$ the \emph{uniform
covering number} of the class $\mathcal{K}$ for $\e>0$, and the
class of norms $\{L_2(\p)\}$, with $\p$ varying in the set of all
probability measures on ${\R}^{d}$ (for more details, see, e.g.,
pp. 83-84 in \cite{vdVW96}). Further set $|{\bf s}|=\max_{1\leq
j\leq d}|s_j|$ for all ${\bf s}\in\R^d$. We will work under the
following assumptions $(K.1$-$2$-$3)$.
\begin{tabbing}
$(K.1)\;\;\;\;$\= $(i)\quad$ \=
$\lim_{|\bf{u}|\tend0}\int_{\R^d}\big(K(\bx)-K(\bx+{\bf
u})\big)^2d\bx=0$.\\[0.2cm]
\>$(ii)$  \>  $\lim_{\lambda\tend1}\int_{\R^d}\big(K(\lambda\bx)
-K(\bx)\big)^2d\bx=0$.\\[0.2cm]
$(K.2)$\>\ $(i)$ \> For some $0<\kappa<\infty$
$K({\bf s})=0$, for $|{\bf s}|\geq \frac{\kappa}{2}$.\\[0.2cm]
\> $(ii)$  \>   For some constant $0<C_K<\infty$, $\sup_{\bx\in\R^d}|K(\bx)|\leq C_K$.\\[0.2cm]
$(K.3)$\>\ $(i)$ \> $\exists\ C>0, v>0, \forall\ 0<\e<1,
\NN(\e,\KK)\leq C\e^{-v}.$\\[0.2cm]
\> $(ii)$  \> $\mathcal{K}$ is pointwise measurable.
\end{tabbing}
\begin{rem}
$(i)\;\;$ It is easily checked that the hypotheses $(K.1)$ and
$(K.3)$ are especially fulfilled if the kernel $K$ can be written
$K(\bx)=\phi(P(\bx))$, where $P$ is a polynomial in $d$ real
variables and $\phi$ a real function with \emph{bounded
variations}.\vskip3pt

\noindent $(ii)\;$ A class $\mathcal{K}$ fulfilling the assumption
$(K.3)(i)$ is said to admit a polynomial uniform covering number.
\end{rem}

\noindent In what follows, we will make use of an auxiliary
function $\{\Theta ({\bx}):{\bx}\in I\}$, assumed to be continuous
and positive on  $I$. Moreover, we will assume that the quantity
$\Theta_{n}({\bx})$ is a consistent estimator of $\Theta(\bx)$, in
the sense that, with probability one,
\begin{tabbing}
$(\Theta.1)\quad$ \=
$\displaystyle\lim_{n\tend\infty}\;\sup_{{\bx}\in
I}\;\Big|{\displaystyle \Theta_n({\bx})\over \displaystyle
\Theta({\bx})}-1\Big|=0$.
\end{tabbing}

\noindent Following the ideas of Einmahl and Mason \cite{EM00}
and Deheuvels and Mason \cite{DM04}, we
will study the uniform convergence to $0$ of
$\widehat{m}^{\star}_{\psi,n,h}$ centered by the following
centering factor.
\begin{equation*}
\widehat{\E}\,m_{\psi;n}({\bx};h) =\E\Big\{\psi ({Y})K\Big(
\frac{{\bx}-\boldsymbol{X}}{h}\Big) \Big\} \Big/
\E\Big\{K\Big(\frac{{\bx}-\boldsymbol{X}}{h}\Big) \Big\}.
\end{equation*}

\begin{rem}
To motivate this choice, we stress the reader attention on the
fact that the remaining bias type term, that is
$\widehat{\E}\,m_{\psi;n}({\bx};h)-m_{\psi}({\bf x})$, can be
neglected under some general regularity assumptions (see Section 2
in \cite{DM04} and Theorem \ref{theo_carbo_pourIC} below)
\end{rem}

\noindent Further introduce for ${\bx}\in I^\alpha$ and
$\psi\in\FF$,
\begin{eqnarray}
\sigma_{\psi }^{2}({\bf x})=\E\{\psi^2 (Y)/[1-G(Y)]\ |\ {\bf
X}={\bf x}\} - m^2_{\psi}({\bf x}). \label{sigma}
\end{eqnarray}

\begin{rem}
Note that, under $(F.1$-$2$-$4)$, the function $\sigma _{\psi }$
(as well as the function $m_\psi$) introduced above is continuous
on $I^{\alp}$ (see, e.g., Section A.3 in \cite{DM04} for a
complete demonstration of this result). It is also noteworthy
that, in the uncensored case, i.e. when $G(y)=0$ for all $y<T_F$,
$\sigma_{\psi }^{2}({\bf x})$ is the conditional variance of
$\psi(Y)$ given ${\bf X}={\bf x}$.
\end{rem}

\section{Main result}\label{section_res_cens}

\setcounter{equation}{0} \setcounter{theo}{0} \setcounter{lem}{0}
\setcounter{cor}{0} \setcounter{rem}{0} \setcounter{prop}{0}
\setcounter{fact}{0}

\noindent We have now all the ingredients to state our main
result, captured in Theorem \ref{theo_carbo} below.

\begin{theo}\label{theo_carbo}
Let $\{h'_n\}_{n\geq 1}$ and $\{h''_n\}_{n\geq 1}$ be two
sequences of positive constants fulfilling the hypotheses
$(H.1$-$2$-$3)$ with $0<h'_n\leq h''_n<1$. Under the assumptions
$({\bf A})$, $(\mathcal{I})$, $(F.1$-$2$-$3$-$4$-$5)$,
$(K.1$-$2$-$3)$, and $(\Theta.1)$, we have with probability one,
\begin{eqnarray}
&&\lim_{\ninf} \suph \suppsi \sup_{{\bx}\in I} \frac{\sqrt{nh^d}
\pm
\Theta_{n}({\bx})\Big\{\widehat{m}^{\star}_{\psi,n,h}({\bx})-\widehat{
\E}\,m_{\psi ;n}({\bx};h)\Big\}}{\sqrt{2\log (1/h^d)}}
 \nonumber
\\[0.08in]
&&= \Bigg\{ \int_{\R^d}K^{2}(\boldsymbol{u})d\boldsymbol{u}\
\sup_{{\bx}\in I} \frac{\Theta^{2}({\bx})\suppsi
\sigma^2_\psi({\bx})}{f_{\boldsymbol{X}}({\bx})} \Bigg\}^{1/2}.
\label{a24a_carbo}
\end{eqnarray}
\end{theo}

\noindent It is noteworthy that in the uncensored case, i.e. when
$G(y)=0$ for all $y<T_F$, the conditions $({\bf A})(ii)$ (and then
the condition $({\bf A})$) , $({\cal I})$ and $(F.3)$ are
automatically fulfilled, in such a way that Theorem
\ref{theo_carbo} reduces in this case to a complement of Corollary
1 in \cite{EM00}, Corollary 3.3 in \cite{DM04} or Theorem 2 in
\cite{EM05}.\vskip5pt

\noindent The proof of Theorem \ref{theo_carbo} is postponed to
Section \ref{section_proofs}. In the following Section
\ref{section_corollaires}, we present some direct consequences of
Theorem \ref{theo_carbo}.

\section{Corollaries-Applications}\label{section_corollaires}

\subsection{Corollaries} \setcounter{equation}{0} \setcounter{theo}{0}
\setcounter{lem}{0} \setcounter{cor}{0} \setcounter{rem}{0}
\setcounter{prop}{0} \setcounter{fact}{0}

In this Section, we show how Theorem \ref{theo_carbo} can be used
$(i)-$ to obtain a uniform law of the logarithm for an estimator
of the conditional distribution function and $(ii)-$ to establish
the uniform-in-bandwidth consistency for some estimates of the
conditional density and conditional hazard rate functions. The
proofs of the corresponding results, captured in Corollaries
\ref{cor_fdr}, \ref{cor_density} and \ref{cor_hazard} below, are
postponed to Section \ref{section_proof_cor}. \vskip3pt

\noindent  Consider the following estimate of the conditional
distribution function $F(t;\bx):= \p(Y\leq t|\bX=\bx)$, along with
the corresponding centering term,
\begin{equation*}
\Fnhhstar(t;{\bx}) := \sum_{i=1}^n W_{n,h,i}({\bx})\frac{\deli
\1_{\{Z_i\leq t\}}}{\Gnstar(Z_i)} \qquad\mbox{and}\qquad
F_{h}(t;{\bx}) := \frac{\E\Big\{\1_{\{Y\leq t\}}K\Big(
{{\frac{{\bx}-\boldsymbol{X}}{h}}}\Big) \Big\}}{\E\Big\{K\Big(
\frac{{\bx}-\boldsymbol{X}}{h}\Big) \Big\}}.
\end{equation*}
\begin{cor} \label{cor_fdr}
Let $\{h'_n\}_{n\geq 1}$ and $\{h''_n\}_{n\geq 1}$ be two
sequences of positive constants fulfilling the hypotheses
$(H.1$-$2$-$3)$ with $0<h'_n\leq h''_n<1$. Under the assumptions
$(\mathcal{I})$, $(F.1$-$2$-$3)$, $(K.1$-$2$-$3)$, and
$(\Theta.1)$, we have, for all $\tau_0<T_H$, with probability one,
\begin{eqnarray}
&&\lim_{\ninf} \suph  \sup_{{\bx}\in I} \frac{\sqrt{nh^d} \pm
\Theta_{n}({\bx})\sup_{t\leq\tau_0}\Big\{\Fnhhstar(t;{\bx})-
F_{h}(t;{\bx})\Big\}}{\sqrt{2\log (1/h^d)}}
 \nonumber
\\[0.08in]
&&= \Bigg\{ \int_{\R^d}K^{2}(\boldsymbol{u})d\boldsymbol{u}\
\sup_{{\bx}\in I} \frac{\Theta^{2}({\bx})\sup_{t\leq\tau_0}
\sigma_{\1_{[0,t]}}^{2}({\bx})}{f_{\boldsymbol{X}}({\bx})}
\Bigg\}^{1/2}. \label{a24b_carbo}
\end{eqnarray}
\end{cor}
\noindent To establish the next corollaries, we will work under
the following additional assumption.
\begin{tabbing}
$(D)\;\; $ \=\  $T_H<\infty$ and the derivatives of order one of
$f_{\bX}$ and $f_{\bX,Y}$
 exist and are \\
\>\  bounded  by a common constant $B_d$ on $I^{\alpha}$ and
$I^{\alpha}\times(-\infty,T_H)$ respectively.
\end{tabbing}
\noindent Moreover, denote by $\ell>0$ an additional bandwidth. As
for the conditional density $f(t;\bx):=f_{\bX, Y}(\bx,
t)/f_{\bX}(\bx)$, we consider the following estimator
\begin{equation*}
\fcond(t;{\bx}) := \sum_{i=1}^n \varpi_{n,h,i}({\bx})\frac{\deli
\1_{\{Z_i\in[t-\frac{\ell}{2};~t+\frac{\ell}{2}]\}}}{\ell\Gnstar(Z_i)}.
\end{equation*}
\begin{cor} \label{cor_density}
Let $\{\ell'_n\}_{n\geq 1}$ and $\{\ell''_n\}_{n\geq 1}$ be two
sequence of positive constants such that $\ell''_n\rightarrow0$,
$\ell''_n\geq \ell'_n\geq\frac{\sqrt{2\log(1/h_n'^d)}}{nh_n'^d} $.
Further assume that the assumptions of Corollary \ref{cor_fdr} and
$(D)$ hold. Then we have, with probability one,
\begin{eqnarray}
\lim_{\ninf} \suph \sup_{\ell\in [\ell'_n,\ell''_n]}
\sup_{{\bx}\in I} \sup_{t\leq\tau_0}\Big\{\fcond(t;{\bx}) -
f(t;\bx)\Big\} =0.
\end{eqnarray}
\end{cor}
\noindent Now, turning our attention to the conditional hazard
rate function $\lambda(t;\bx):=f(t;\bx)/[1-F(t;\bx)]$, we
introduce the following estimator
\begin{equation*}
\widehat{\lambda}_{n,h,\ell}^{\star}(t;{\bx}):=\frac{\widehat{f}_{n,h,\ell}^{\star}(t;{\bx})}{1-\Fnhhstar(t;{\bx})}.
\end{equation*}
\begin{cor}\label{cor_hazard}
Under the assumptions of Corollary \ref{cor_density}, we have,
with probability one,
\begin{eqnarray}
&&\lim_{\ninf} \suph \sup_{\ell\in [\ell'_n,\ell''_n]}
\sup_{{\bx}\in I}
\sup_{t\leq\tau_0}\Big\{\widehat{\lambda}_{n,h,\ell}^{\star}(t;{\bx})
-\lambda(t;\bx)\Big\} =0).
\end{eqnarray}
\end{cor}

\subsection{Almost sure asymptotic simultaneous
confidence bands for the true regression function}
\label{section_IC}

\subsubsection{A necessary result}
In this section, we first establish a direct consequence of
Theorem \ref{theo_carbo} which will enable us to build confidence
bands for the theoretical regression function. The following
notations and hypotheses will be needed.\vskip3pt

\noindent In this paragraph, we will impose the following
regularity conditions upon the functions $f_{\boldsymbol{X}}$ and
$f_{\boldsymbol{X},Y}$.
\begin{tabbing}
$(F.6)$ \=$(i)\quad$ \= $f_{\boldsymbol{X}}$ is three times continuously differentiable on  $I^\alpha$.\\[0.2cm]
 \>$(ii)$ \> $f_{\boldsymbol{X},Y}$ is three times continuously differentiable on  $I^\alpha\times \R$.
\end{tabbing}
To ensure that the bias-type term may be neglected, we will work
in this section with kernels satisfying the assumption $(K.4)$
below.
\begin{tabbing} $(K.4)\;\;$\=\
$\int_{\R^d}u_{1}^{j_1}\ldots u_{d}^{j_d}K({\bf u}) d{\bf
u}=0,\quad j_1,\ldots,j_d\geq 0,\quad j_1+\ldots+j_d=0,1,2$.
%\>\ $(ii)$\>\ $\int_{\R^d}|u_{1}^{j_1}\ldots u_{d}^{j_d}|K({\bf
%u}) d{\bf u}<\infty,\quad j_1,\ldots,j_d\geq 0,\quad
%j_1+\ldots+j_d=3$.
\end{tabbing}
\noindent Further denote by ${\cal V}_I$ the Lebesgue measure of
$I$ and set
\[\log_{\theta,K}(u):=\log\Lpar\theta\vee
u\Lac\displaystyle\int_{\R^d}K^2(\boldsymbol{t})d\boldsymbol{t}\Rac\Rpar,\]
where $\theta>1$ is a fixed constant. We refer to Remark 1.4 in
\cite{DM04} for discussions about the introduction of these
quantities along with some relevant choices for $\theta$.

\noindent Finally, given a sequence $\{h_n\}_{n\geq 1}$ of
positive constants fulfilling the conditions $(H.1$-$2$-$3)$,
consider a sequence of possibly data-driven bandwidths $H_n({\bf
x})$ such that the assumption $(B.1)$ below holds, for two given
constants $0<c_1\leq c_2<\infty$.
\begin{tabbing}
$(B.1)\;\;$ \= $\displaystyle c_1h_n\leq \inf_{\boldsymbol{x}\in
I}H_n(\boldsymbol{x})\leq \sup_{\boldsymbol{x}\in
I}H_n(\boldsymbol{x})\leq c_2h_n$\quad almost surely as $\ninf$.
\end{tabbing}

\begin{theo}\label{theo_carbo_pourIC}
Given two constants $0<A<\infty$ and $1/(4+d)\leq\del_0<1$, assume
that $H_n(\bx)$ is a bandwidth function fulfilling the assumption
$(B.1)$ with $h_n=An^{-\del_0}$ . Then, under the assumptions
$({\bf A})$, $(\mathcal{I})$, $(F.1$-$2$-$3$-$4$-$5$-$6)$,
$(K.1$-$2$-$3$-$4)$ and $(\Theta.1)$, we have, with probability
one,
\begin{eqnarray}
&&\lim_{\ninf} \suppsi\sup_{\boldsymbol{x}\in I} \frac{\sqrt{n
H^d_n(\boldsymbol{x})} \pm
\Theta_{n}(\boldsymbol{x})\Big\{\widehat{m}^{\star}_{\psi,n,
H_n(\boldsymbol{x})}(\boldsymbol{x}) -
m_{\psi}(\boldsymbol{x})\Big\}}{\sqrt{2\log_{\theta,K} ({\cal
V}_I/H^d_n(\boldsymbol{x}))}}
 \nonumber
\\[0.08in]
&&= \Bigg\{ \int_{\R^d}K^{2}(\boldsymbol{t})d\boldsymbol{t}\
\sup_{\boldsymbol{x}\in I}
\frac{\Theta^{2}(\boldsymbol{x})\suppsi\sigma_{\psi
}^{2}(\boldsymbol{x})}{f_{\boldsymbol{X}}(\boldsymbol{x})}
\Bigg\}^{1/2}. \label{a24am1}
\end{eqnarray}
\end{theo}

\noindent The proof of Theorem \ref{theo_carbo_pourIC} follows
from Theorem \ref{theo_carbo} along the same lines as the first
part of Theorem 1.1 in \cite{DM04} is shown to be a consequence of
their Corollary 3.3. We omit the details of this book-keeping
argument.

\begin{rem}
The assumption $(B.1)$ above is the almost sure version of
condition $(B.1)$ in \cite{DM04}. Theorem \ref{theo_carbo} should
allow to treat even more general random bandwidths, especially
bandwidths $H_n$ such that $h'_n\leq \inf_{\boldsymbol{x}\in
I}H_n(\boldsymbol{x})\leq \sup_{\boldsymbol{x}\in
I}H_n(\boldsymbol{x})\leq h''_n$ almost surely as $\ninf$,
whenever $\{h'_n\}_{n\geq 1}$ and $\{h''_n\}_{n\geq 1}$ are two
sequences of positive constants fulfilling $(H.1$-$2$-$3)$.
However, some regularity conditions have to be imposed upon the
function $H_n(\cdot)$, otherwise the corresponding version of
Theorem \ref{theo_carbo_pourIC} can not be directly derived from
Theorem \ref{theo_carbo} (the equality in (\ref{a24am1}) becoming
an inequality). Alternatively, this generalization could be
derived from some functional limit law associated to Theorem
\ref{theo_carbo}. We especially refer the reader interested in
functional limit laws of the logarithm to \cite{DE00AP} or
\cite{Mason04} and the relevant references therein.
\end{rem}

\subsubsection{Construction of the bands}
Here we give an example illustrating how Theorem
\ref{theo_carbo_pourIC} may be used to construct simultaneous
confidence bands for $m_\psi$. Consider a random sequence of
functions $L_n(\bx)$ such that for all $0<\e<1$, there exists
almost surely an $n_0=n_0(\e)$ such that, for all $n\geq n_0$,
uniformly over $\boldsymbol{x}\in I$,
\begin{equation}
\begin{split}
&m_{\psi}(\boldsymbol{x})\in [\widehat{m}^{\star}_{\psi,n,
H_n(\boldsymbol{x})}(\boldsymbol{x}) \pm
(1+\e)L_n(\boldsymbol{x})],  \\
& m_{\psi}(\boldsymbol{x})\notin [\widehat{m}^{\star}_{\psi,n,
H_n(\boldsymbol{x})}(\boldsymbol{x}) \pm
(1-\e)L_n(\boldsymbol{x})]. \label{def_IC1}
\end{split}
\end{equation}
Whenever $(\ref{def_IC1})$ is fulfilled for all $0<\e<1$, we will
say that the intervals
\begin{equation}
[A_n(\boldsymbol{x}), B_n(\boldsymbol{x})] =
[\widehat{m}^{\star}_{\psi,n, H_n(\boldsymbol{x})}(\boldsymbol{x})
- L_n(\boldsymbol{x}), \widehat{m}^{\star}_{\psi,n,
H_n(\boldsymbol{x})}(\boldsymbol{x}) + L_n(\boldsymbol{x})]
\end{equation}
provide asymptotic simultaneous confidence bands for $m_\psi(\bx)$
over $\bx\in I$.\vskip5pt

\noindent It is easy to check that the quantity
\begin{eqnarray}
{\sigma}_{\psi;n}^{\star 2}(\boldsymbol{x};H_n(\boldsymbol{x})) =
\sum_{i=1}^{n}\frac{\deli (\psi(Z_{i}))^{2}}{(1-\Gnstar(Z_i))^2}
\varpi_{n,H_n(\boldsymbol{x}),i}(\boldsymbol{x}) -
\widehat{m}^{\star 2}_{\psi,n,H_n(\boldsymbol{x})}(\boldsymbol{x})
\label{sigmacens}
\end{eqnarray}
is a consistent estimate of $\sigma_\psi(\bx)$, uniformly over
$\psi\in\FF$ and $\bx\in I$. Therefore, setting
\begin{equation}
L_n(\boldsymbol{x})= \left\{ \frac{2\log_{\theta,K}({\cal V}_I /
H^d_n(\boldsymbol{x}) )}{n
H^d_n(\boldsymbol{x})}\times\frac{\widetilde{\sigma} _{\psi;n
}^{\star
2}(\boldsymbol{x},H_n(\boldsymbol{x}))}{f_{\boldsymbol{X};n}(\boldsymbol{x},H_n(\boldsymbol{x}))}\right\}^{1/2}
\times\left[\int_{\R^d}K^2(\boldsymbol{t})d\boldsymbol{t}\right]^{1/2},
\label{ICintermed3}
\end{equation}
with
\begin{equation*}
f_{{\bf X};n}({\boldsymbol x};H_n(\boldsymbol{x}))
={\frac{1}{nH^d_n(\boldsymbol{x})}}\sum_{i=1}^{n}K\Big(\frac{{\boldsymbol
x}-{\bf X}_{i}}{H_n(\boldsymbol{x})} \Big),
\end{equation*}
it is easily derived from Theorem \ref{theo_carbo_pourIC} that the
intervals
\begin{equation}
[\widehat{m}^{\star}_{\psi,n, H_n(\boldsymbol{x})}(\boldsymbol{x})
- L_n(\boldsymbol{x}), \widehat{m}^{\star}_{\psi,n,
H_n(\boldsymbol{x})}(\boldsymbol{x}) + L_n(\boldsymbol{x})]
\label{ICfin}
\end{equation}
provide asymptotic simultaneous confidence bands for $m_\psi(\bx)$
over $\bx\in I$, in the sense made precise above.

\subsubsection{Illustration : a simple simulation study} In this
paragraph, the confidence bands introduced above are constructed
on simulated data. We worked with a sample size $n=2000$, and
considered the case where ${\bf X}=X\in\R$ (i.e. $d=1$) was such
that $X\sim\mathcal{N}(0,1)$, where $\mathcal{N}(0,1)$ stands for
the gaussian distribution with mean 0 and standard deviation 1.
Set $p(x)=0.25 + 0.5\times\cos^2(x)$. We selected
$\psi=\1_{\{.\leq 0.9\}}$, and considered the model
$\E[\psi(Y)|X=x]= p(x)$. Under this model, the variable $Y$ was
simulated as follows. For each integer $1\leq i \leq n$, let
$p_i=p(x_{i})$ where $x_{i}$ is the observed value of the variable
$X_i$. Note that $0< p_i< 1$ for every $1\leq i \leq n$. Each
$Y_i$ was then generated as one $\mathcal{U}(0.9-p_i, 1+0.9-p_i)$
variable, where $\mathcal{U}(a,b)$ stands for the uniform
distribution on (a,b). Following this proceed ensured that
$\p(Y_i\leq 0.9|X_i=x_i)=p_i=p(x_{i})$. Regarding the censoring
variable, we generated an i.i.d. sample $C_1, ..., C_n$ such that
$C_i \sim\mathcal{U}(0,1)$. This choice yielded, \emph{a
posteriori}, $\p(\delta=1)\simeq 0.2$. As for $K$, we opted for
the Epanechnikov kernel. Moreover, we selected $H_n(\bx)=h=0.15$
and $H_n(\bx)=h=0.20$.
\begin{figure}
\subfigure[h=0.15]{
\includegraphics[width=0.45\textwidth]{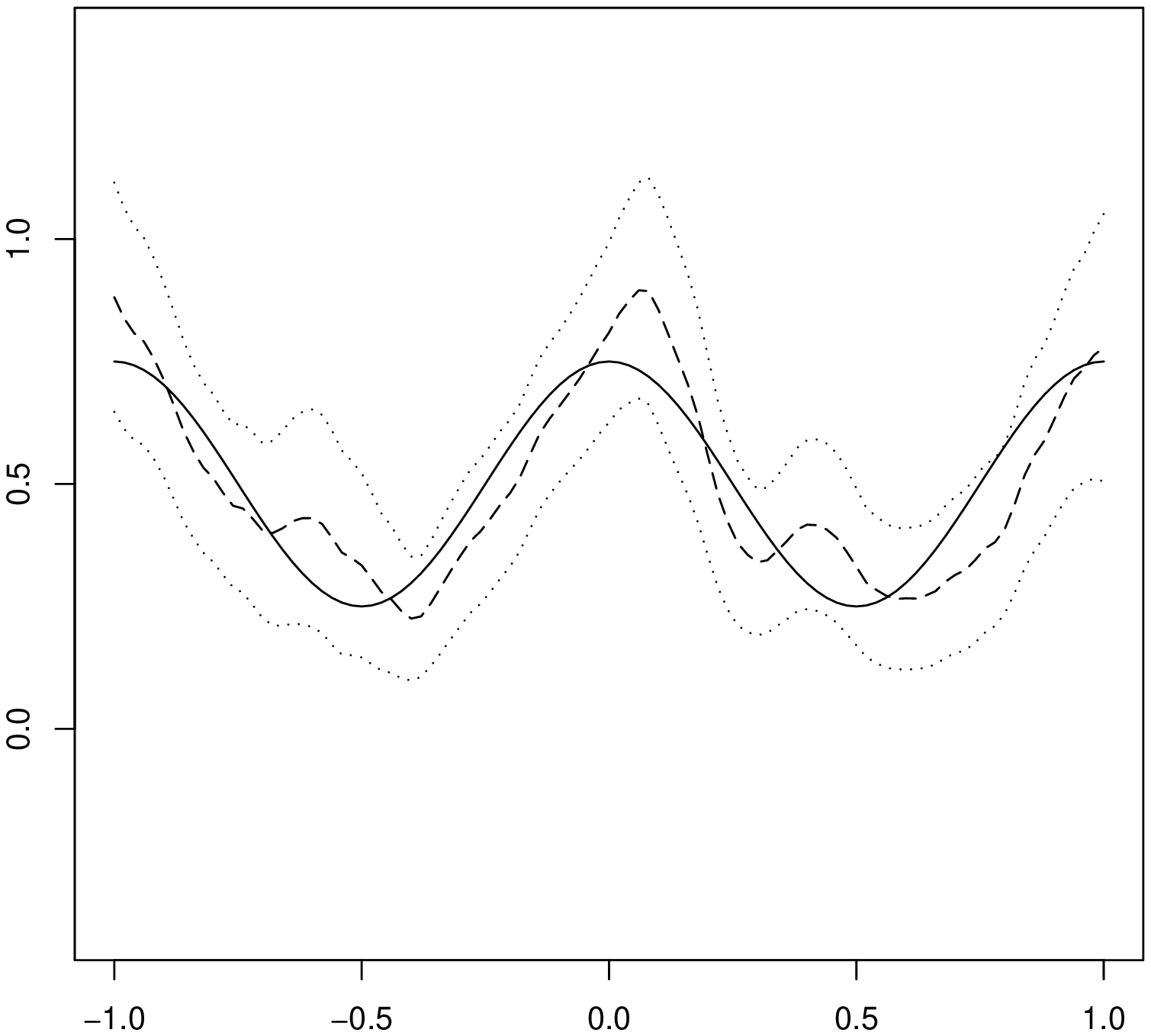}}
\hspace{.3in} \subfigure[h=0.20]{
\includegraphics[width=0.45\textwidth]{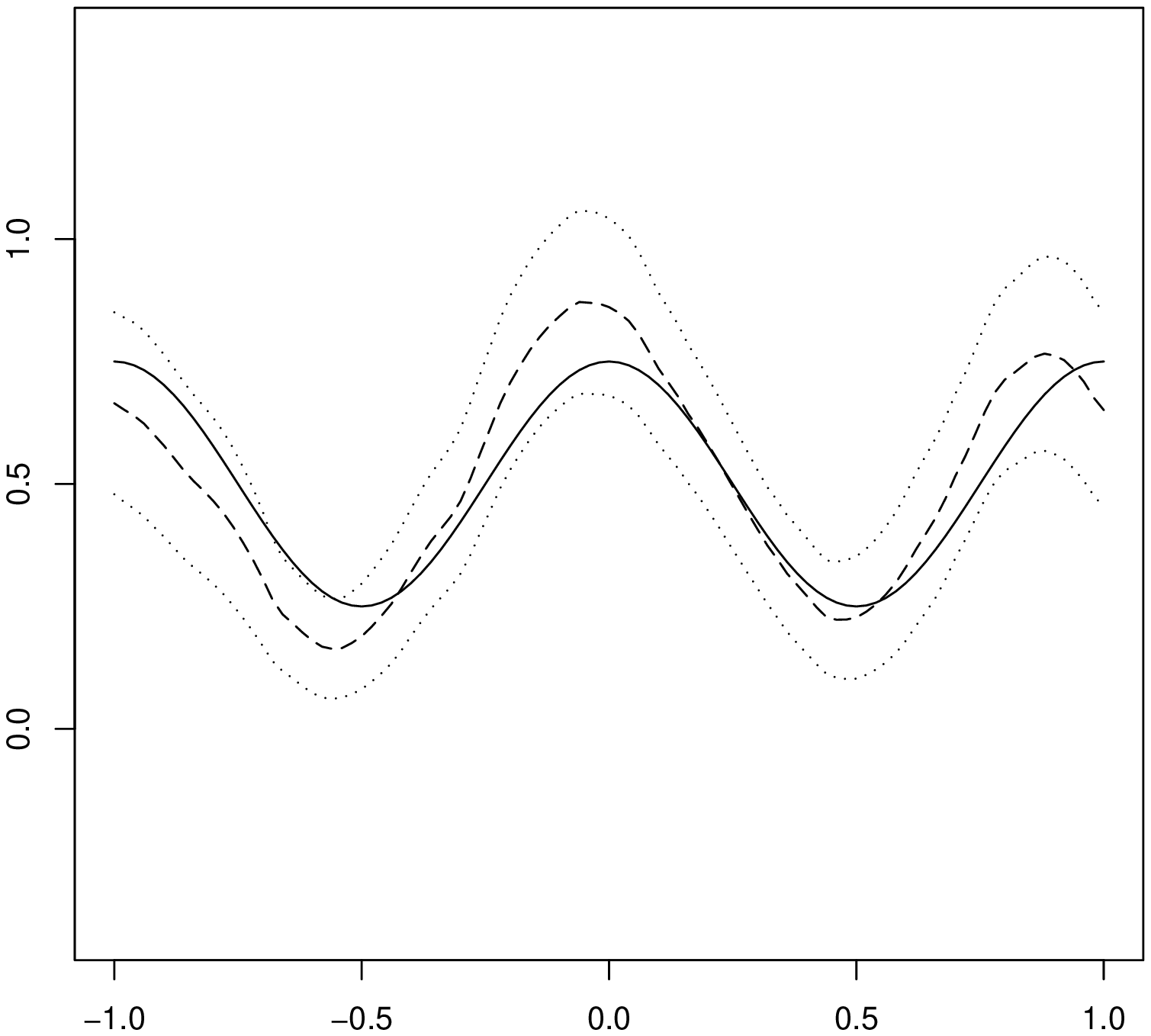}}
\caption{Results of the simulation study for (a) $h=0.15$ and (b)
$h=0.20$ : true additive components (solid line), their
estimations (dashed line), and the associated confidence bands
(dotted line)}\label{graph}
\end{figure}
\noindent Results are presented in Figure \ref{graph}. The
confidence bands appeared to be adequate, in the sense that they
contained the true value of the regression function for every
$x\in[-1,1]$. The fact that the true function did not belong to
our bands for some points was expected : it is due to the
$\varepsilon$ term in (\ref{def_IC1}).\vskip5pt

\begin{figure}
\subfigure[h=0.15]{
\includegraphics[width=0.45\textwidth]{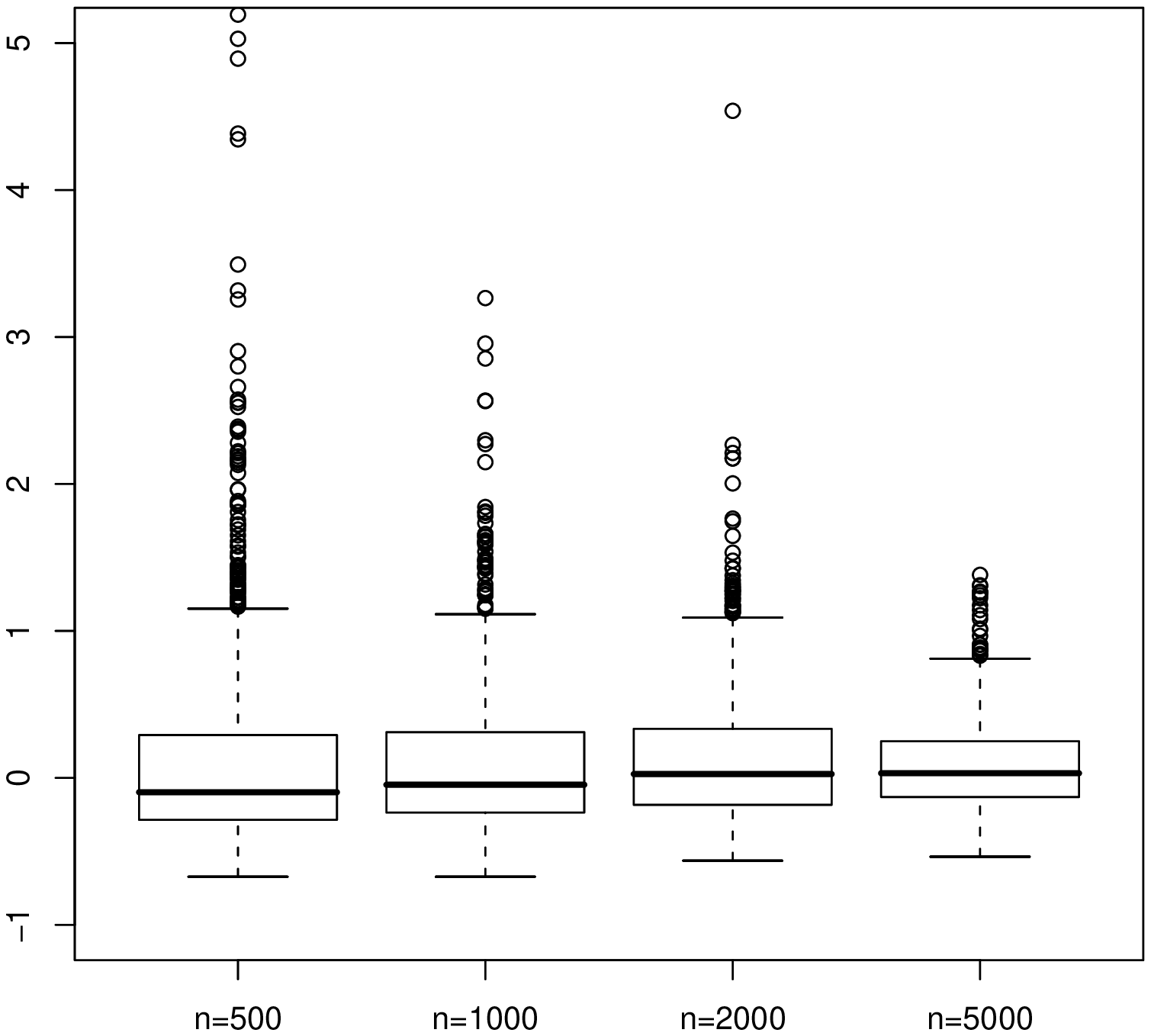}}
\hspace{.3in}\subfigure[h=0.20]{
\includegraphics[width=0.45\textwidth]{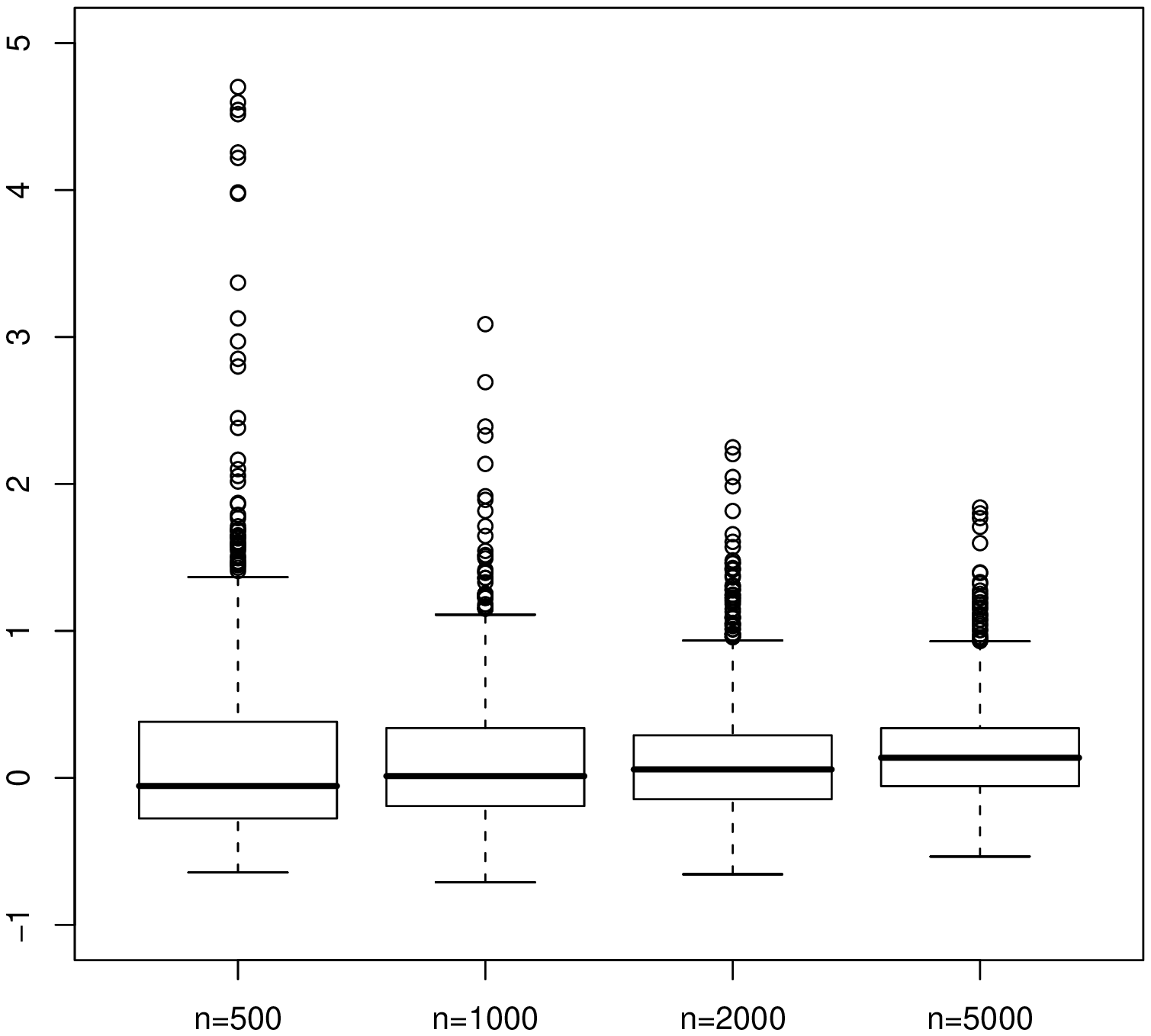}}
\caption{Boxplots of the $\e_1(h,n)$ distribution for various
sample sizes $n$ and for (a) $h=0.15$ and (b)
$h=0.20$.}\label{graph2}
\end{figure}

\noindent We also performed a simulation study in order to
empirically evaluate the finite sample behavior of the
distribution of the $\e$ term involved in (\ref{def_IC1}). Given a
sample of size $n$ (simulated in the same way as above), and an
estimate $\widehat{m}^{\star}_{\psi,n, h}({x})$ built on this
sample, consider the quantity
\[ \e_1(h,n) = |\widehat{m}^{\star}_{\psi,n,h}({x}_0)-m_{\psi}({x}_0)| - L_n({x}_0),\]
where $L_n(x)$ is defined as in (\ref{ICintermed3}) and
\[x_0=\mbox{argmax }\{|\widehat{m}^{\star}_{\psi,n,
h}({x})-m_{\psi}({x})|,\ x\in [-1,1]\}.\] Then, by simulating 2000
samples of size $n$, we were able to estimate the distribution of
$\e_1(h,n)$, for various sizes $n$ and bandwidths $h$. Our results
are presented in Figure \ref{graph2}. As expected, the
distribution of $\e_1(h,n)$ is all the more concentrated around 0
as $n$ is high. Moreover, the bandwidth choice is crucial since
$\e_1(h,n)$ appears not to be centered around 0 for some pairs
$(h,n)$. This highlights the need to build some procedures aiming
at optimally selecting the bandwidth (in a sense to be made
precise). This particular problem will be addressed elsewhere. It
is also noteworthy that a more formal description of the
asymptotic behavior of the $\e$ term involved in (\ref{def_IC1})
could be achieved by studying the rate of coverage pertaining to
Theorem \ref{theo_carbo} and Theorem \ref{theo_carbo_pourIC} (see
\cite{Berthet1} and \cite{Berthet2} for examples of such results).

\section{Proofs}\label{section_proofs}

\setcounter{equation}{0} \setcounter{theo}{0} \setcounter{lem}{0}
\setcounter{cor}{0} \setcounter{rem}{0} \setcounter{prop}{0}
\setcounter{fact}{0}

To prove our results, we will first establish a general result for
an estimator of the regression function in the uncensored case.
Then, we will make use of the relation $(\ref{key_idea})$ to treat
the censored case when the function $G$ is known. Finally, strong
consistency results for the Kaplan-Meier estimator will be
employed to cope with the general censored case.

\subsection{The uncensored case}
\noindent Let $(\boldsymbol{X}_{1}, {\bf Y}_{1}),
(\boldsymbol{X}_{2}, {\bf Y}_{2}), \ldots$, be a sequence of
independent and identically distributed [i.i.d.] replica of the
random pair $(\boldsymbol{X},{\bf Y}) \in{\R}^d \times {\R^q},\,
d, q \geq 1$. Consider a pointwise measurable class $\FF_q$ of
real measurable functions defined on $\R^q$. We will assume that
$\FF_q$ forms a \emph{VC subgraph class} and consider the
conditional expectation of $\Psi({\bf Y})$ given
$\boldsymbol{X}={\bx}$, for $\Psi\in\FF_q$,
\begin{equation}
m_{\Psi}({\bx})={\E}\big(\Psi ({\bf Y})\mid
\boldsymbol{X}={\bx}\big). \label{a1}
\end{equation}
\noindent In this uncensored setting, we will especially work
under the assumptions $(F.I$-$IV$-$V)$ below.
\begin{tabbing}
$(F.I)\;\;$ \=\ For all ${\bx}\in I^{\alp}$,
$\displaystyle\lim_{{\bx}'\rightarrow {\bx};{\bx}'\in
I^{\alp}}f_{\boldsymbol{X},{\bf Y}}({\bx}',{\bf
y})=f_{\boldsymbol{X},{\bf Y}}({\bx},{\bf y})$ for almost every
${\bf y}\in\R^q$. \\
$(F.IV)$ \>\ The class of functions $\FF_q$ is bounded.\\
$(F.V)$ \>\ The class of functions $\MM_q:=\{m_{\Psi
}/f_{\boldsymbol{X}}, \Psi\in\FF_q\}$ is relatively compact with
\\ \>\ respect  to the sup-norm  topology on $I^\alpha$.
\end{tabbing}
Under $(F.I$-$2$-$IV)$, the conditional variance $\sigma
_{\Psi}^{2}({\bx})$ of $\Psi ({\bf Y})$ given
$\boldsymbol{X}={\bx}$ is defined, for $\Psi\in\FF_q$, by
\begin{equation}
\widetilde{\sigma}_{\Psi }^{2}({\bx})={\rm Var}(\Psi({\bf Y})\mid
\boldsymbol{X}={\bx})={\frac{1}{f_{\boldsymbol{X}}({\bx})}}\int_{\R^q}\Big(\Psi(\bf
y)-m_{\Psi }({\bx})\Big)^{2}f_{\boldsymbol{X},{\bf Y}}({\bx},{\bf
y})d{\bf y}. \label{a5}
\end{equation}
\noindent Introduce the kernel estimator of $m_{\Psi}({\bx})$
defined on $I$, for any $\inth$, by
\begin{equation*}
m_{\Psi,n,h}({\bx}) := \sum_{i=1}^n \Psi({\bf Y}_i)
\varpi_{n,h,i}({\bx}),
\end{equation*}
along with the following centering term
\begin{equation*}
\widehat{\E}\,m_{\Psi;n}({\bx};h) =\E\Big\{\Psi (\bf{Y})K\Big(
{{\frac{{\bx}-\boldsymbol{X}}{h}}}\Big) \Big\} \Big/
\E\Big\{K\Big( {{\frac{{\bx}-\boldsymbol{X}}{h}}}\Big) \Big\}.
\end{equation*}

\subsubsection{Results} \noindent We have now all the ingredients
to state the result corresponding to the uncensored case. As
mentioned above, Theorem \ref{theo_reg} below is a complement of
the successive results of Einmahl and Mason \cite{EM00}, Deheuvels
and Mason \cite{DM04} and Einmahl and Mason \cite{EM05}.
\begin{theo}\label{theo_reg}
Let $\{h'_n\}_{n\geq1}$ and $\{h''_n\}_{n\geq1}$ be two sequences
of positive constants fulfilling the conditions $(H.1$-$2$-$3)$,
with $0<h'_n\leq h''_n<1$. Under the hypotheses
$(F.I$-$2$-$IV$-$V)$ and $(K.1$-$2$-$3)$, we have with probability
one,
\begin{eqnarray}
&&\lim_{\ninf} \suph \Big\{\frac{nh^d}{2\log (1/h^d)}\Big\}^{1/2}
\sup_{\Psi\in \FF_q} \sup_{{\bx}\in
I}\pm\Theta_n({\bx})\Big\{m_{\Psi;n}({\bx};h)-\widehat{
\E}\,m_{\Psi ;n}({\bx};h)\Big\}
 \nonumber
\\
&&\hskip30pt= \Bigg\{ \int_{\R^d}K^{2}({\bf t})d{\bf t}\
\sup_{{\bx}\in I} \frac{\Theta(\bx)\sup_{\Psi\in \FF_q}
\widetilde{\sigma}_{\Psi }^{2}({\bx})}{f_{\boldsymbol{X}}({\bx})}
\Bigg\}^{1/2}. \label{a24a}
\end{eqnarray}
\end{theo}

\noindent In Section \ref{section_final_proof}, Theorem
\ref{theo_reg} will be shown to be a direct consequence of the
technical result captured in Theorem \ref{theo_W} in Section
\ref{vvvv1}. In the following Section \ref{section_proof_carbo},
we first show that Theorem \ref{theo_carbo} naturally follows from
Theorem \ref{theo_reg}.

\subsection{Proof of Theorem
\ref{theo_carbo}}\label{section_proof_carbo}

\noindent In Proposition \ref{propanalog} below, we first
establish the version of Theorem \ref{theo_carbo} corresponding to
the case where $G$ is known (i.e. with
$\widehat{m}^{\star}_{\psi,n,h}$ replaced by
$\widehat{m}_{\psi,n,h}$). To complete the proof of Theorem
$\ref{theo_carbo}$ in the general case, the consistency of the
Kaplan-Meier estimate will be helpful (see Lemmas \ref{lemApprox1}
and \ref{lemApprox2} below).

\begin{prop}\label{propanalog}
Let $\{h'_n\}_{n\geq 1}$ and $\{h''_n\}_{n\geq 1}$ be two
sequences of positive constants fulfilling $(H.1$-$2$-$3)$, with
$0<h'_n\leq h''_n<1$. Under the hypotheses $(\bf{A})$,
$(\mathcal{I})$, $(F.1$-$2$-$3$-$4$-$5)$, $(K.1$-$2$-$3$) and
$(\Theta.1)$, we have almost surely
\begin{eqnarray}
&&\lim_{\ninf} \suph \suppsi\sup_{{\bx}\in I} \frac{\sqrt{nh^d}
\pm \Theta_{n}({\bx})\Big\{\widehat{m}_{\psi,n,h}-\widehat{
\E}\,m_{\psi ;n}({\bx};h)\Big\}}{\sqrt{2\log (1/h^d)}} \nonumber
\\[0.08in]
&&= \suppsi\Bigg\{
\int_{\R^d}K^{2}(\boldsymbol{t})d\boldsymbol{t}\ \sup_{{\bx}\in I}
\frac{\Theta^{2}({\bx})\widetilde{\sigma}_{\psi
}^{2}({\bx})}{f_{\boldsymbol{X}}({\bx})} \Bigg\}^{1/2}.
\label{analogchap4}
\end{eqnarray}
\end{prop}
\demo\ Recalling the definition (\ref{def_phipsi}) of $\Phi_\psi$,
this function is uniformly bounded under $({\bf A})$ and $(F.4)$,
for all $\psi\in\FF$. This property, when combined with the VC
property of $\FF$, ensures that the class of functions
$\FF_\Phi:=\{\Phi_\psi, \psi\in\FF\}$ has a polynomial uniform
covering number. Similarly, it can be shown that $\FF_\Phi$ is a
pointwise measurable class of functions. Moreover, by $(F.IV)$ and
$(\ref{key_idea})$, the class $\MM'=\{\E(\Phi_\psi(Y,C)|{\bf
X}=\cdot), \psi\in\FF\}$ is almost surely relatively compact with
respect to the sup-norm topology. Therefore, we can apply Theorem
$\ref{theo_reg}$ with ${\bf Y} = (Y,C)$ and $\Psi=\Phi_\psi$. In
this setting, note that $(F.I)$ follows from $({\cal I})$, $(F.1)$
and $(F.3)$. Now, observing that under the assumption $({\cal
I})$, for $j=1,2$,
\begin{eqnarray*}
\E\Lac\Lpar\frac{\delta\psi(Z)} {G(Z)}\Rpar^j\Big|{\bf X}\Rac
%&=&\ \E\Lac\Lpar\frac{\psi(Y)}{1-G(Y)}\Rpar^j \E\big[\1_{\{Y\leq
%C\}}|{\bf X}, Y\big]\Big|{\bf X}\Rac \nonumber \\
=\ \E\lac\frac{(\psi(Y))^j}{\big(G(Y)\big)^{j-1}} \Big|{\bf
X}\rac,
\end{eqnarray*}
the result of Proposition \ref{propanalog} is straightforward.
\findem \vskip5pt

\noindent Lemma \ref{lemApprox1} below enables to complete to
proof of Theorem \ref{theo_carbo} in the case where $({\bf A})(i)$
holds.
\begin{lem}\label{lemApprox1}
Assume the hypotheses of Theorem $\ref{theo_carbo}$ hold with
$({\bf A})(i)$. Then, we have, with probability one,
\begin{equation}
\suph \suppsi\sup_{{\bx}\in I}
|\widehat{m}^{\star(1)}_{\psi,n,h}({\bx}) -
\widehat{m}^{\star}_{\psi,n,h}({\bx})| = o\lpar \sqrt{\frac{\log
(1/h_n)}{nh^d_n}}\rpar, \mbox{ as } \ninf. \label{approx1}
\end{equation}
\end{lem}
\demo\ Keep in mind the definition \ref{poids_carbo} of the
functions $\varpi_{n,h,i}$. Setting $\varpi_{n,I}=\suph\supx
\sum_{i=1}^n |\varpi_{n,h,i}({\bx})|$, it is easy to check that
$\varpi_{n,I}<\infty$. Moreover, observe that, under $({\bf
A})(i)$,
\begin{equation*}
\suppsi \suphx |\widehat{m}^{\star(1)}_{\psi,n,h}({\bx}) -
\widehat{m}^{\star}_{\psi,n,h}({\bx})| \leq  \varpi_{n,I}
\frac{\sup_{t \leq \tau_0}\{\suppsi \psi(t)\times
|\Gnstar(t)-G(t)|\}}{\Gnstar(\tau_0)G(\tau_0)}.
\end{equation*}
Since $\tau_0<T_H\leq T_G$, the law of the logarithm for $\Gnstar$
established in \cite{foldesrejto} ensures that
\[\sup_{t \leq \tau_0} |\Gnstar(t)-G(t)|=O\lpar
\sqrt{\displaystyle \frac{\log_2n}{n}}\rpar\quad \mbox{almost
surely as $\ninf$},\]  which allows to conclude to (\ref{approx1})
under $(H.3)$. \findem \vskip5pt

\noindent Lemma \ref{lemApprox2} below enables to complete to
proof of Theorem \ref{theo_carbo} in the case where $({\bf
A})(ii)$ holds.
\begin{lem}\label{lemApprox2}
Assume the hypotheses of Theorem $\ref{theo_carbo}$ hold with
$({\bf A})(ii)$. Then, we have, with probability one,
\begin{equation}
\suph \suppsi\sup_{{\bx}\in I}
|\widehat{m}^{\star(1)}_{\psi,n,h}({\bx}) -
\widehat{m}^{\star}_{\psi,n,h}({\bx})| = o\lpar \sqrt{\frac{\log
(1/h_n)}{nh^d_n}}\rpar, \mbox{ as } \ninf. \label{approx2}
\end{equation}
\end{lem}
\demo\ Observe that, under $({\bf A})(ii)$,
\begin{equation*}
\suppsi \suphx |\widehat{m}^{\star(1)}_{\psi,n,h}({\bx}) -
\widehat{m}^{\star}_{\psi,n,h}({\bx})| \leq  \varpi_{n,I} \sup_{t
\leq \mathcal{Z}_n} \frac{\suppsi
\psi(t)}{\Gnstar(t)G(t)}\times\sup_{t \leq T_H} |\Gnstar(t)-G(t)|,
\end{equation*}
where $\varpi_{n,I}$ is defined as in the proof of Lemma
\ref{lemApprox1} and $\mathcal{Z}_n:=\max\{Z_i:\delta_i=1, 1\leq
i\leq n\}$. Clearly, $\min(G(\mathcal{Z}_n),
\Gnstar(\mathcal{Z}_n))>0$ for all $n\geq 1$. \vskip3pt

\noindent Suppose that $0<p<1/2$ in $({\bf A})(ii)$. Then, Theorem
2.1 in \cite{ChenLo} ensures that
\[\sup_{t \leq T_H} |\Gnstar(t)-G(t)|=o\big( n^{-p}\big)\quad\mbox{almost surely
as $\ninf$},\] which allows to conclude to (\ref{approx2}) under
$(H.1$-$2$-$3)$ and $({\bf A})(ii)(c)$.\vskip3pt

\noindent The proof in the case where $({\bf A})(ii)$ holds with
$p=1/2$ follows from the same lines, replacing the result of
Chen and Lo \cite{ChenLo} by that of Gu and Lai \cite{GuLai}. Details are omitted.
\findem

\subsection{A useful technical result}\label{vvvv1}

\noindent In Section \ref{section_final_proof} (see, especially,
Lemma \ref{lem_approx_estim_proc_emp_index}), we will show that
Theorem \ref{theo_reg} is a direct consequence  of a general
empirical process result, captured in Theorem \ref{theo_W} below.
This latter result describes the oscillations of a version of the
multivariate empirical process indexed by an appropriate class of
functions. This process is defined in
(\ref{proc_emp_multiva_index_fonctions}) below. For any function
$\Psi\in\FF_q$ and any pair of continuous functions, $c_\Psi$ and
$d_\Psi$, defined on the compact $I\subset\R^d$, set, for all
${\bx}\in I$ and any $h>0$,
\begin{eqnarray}
W_{n,h}({\bx},\Psi ) &=&\sum_{j=1}^{n}\Big(c_\Psi({\bx})\Psi({\bf
Y}_{j})+d_\Psi({\bx})\Big)K\Big(\frac{{\bx}-\boldsymbol{X}_{j}}{h}\Big)
\nonumber \\[0.08in]
&&-n{{\E}}\Big\{ \Big(c_\Psi({\bx})\Psi({\bf
Y})+d_\Psi({\bx})\Big)K \Big(\frac{{\bx}-\boldsymbol{X}}{h}\Big)
\Big\} . \label{proc_emp_multiva_index_fonctions}
\end{eqnarray}
Further denote by $\alpha_n$ the empirical process based on the
observations $(\boldsymbol{X}_1,{\bf Y}_1),$$\ldots,$
$(\boldsymbol{X}_n,{\bf Y}_n)$, and indexed by a class of
functions ${\cal G}$. Namely, for $g \in {\cal G}$, $\alpha_n(g)$
is defined by
\begin{equation}
\alpha_{n}(g)=\frac{1}{\sqrt{n}}\sum_{i=1}^{n}\big(
g(\boldsymbol{X}_{i},{\bf Y}_{i})-{\E}\,g(\boldsymbol{X}_{i},{\bf
Y}_{i})\big). \label{mas1}
\end{equation}
It is noteworthy that, setting for $\Psi\in\FF_q$, ${\bx}\in I$
and $h>0$,
\begin{equation}
\eta_{{\bx},h,\Psi}({\bf u},{\bf v})=\Big(c_\Psi({\bx})\Psi({\bf
v})+d_\Psi({\bx})\Big)
 K\Big(\frac{{\bx}-{\bf
u}}{h}\Big) ,\quad\hbox{for}\quad  {\bf u},{\bf v}\in {\R}^d
\times {\R^q}. \label{hz}
\end{equation}
 the following relation holds,
\begin{equation}
W_{n, h}({\bx},\Psi )=n^{1/2}\alpha_{n}(\eta_{{\bx},h, \Psi}).
\label{wr}
\end{equation}
In view of (\ref{wr}), the introduction of the process $\alpha_n$
provides a suitable and general set-up to study various types of
kernel estimators (especially the density and regression functions
estimators). See \cite{EM00}, \cite{DM04}, \cite{EM05} and
\cite{DonyEM} for details.\vskip5pt

\noindent  For future use, consider the following class of
functions,
\begin{eqnarray}
\mathcal{G}^{'} := \{ \eta_{{\bx}_1,h_1,\Psi_1} -
\eta_{{\bx}_2,h_2,\Psi_2}:~  {\bx}_1, {\bx}_2\in I,~\Psi_1, \Psi_2
\in \mathcal{F}_q, h_1,h_2>0 \}.\label{Gkprime}
\end{eqnarray}
By $(K.3)$, arguing as in pages 17 and 18  of \cite{EM00}, it can
be shown that $\mathcal{G}^{'}$ is a pointwise measurable class of
measurable functions admitting a bounded envelope function and a
polynomial uniform covering number. \vskip5pt

\noindent In other respect, introduce the following classes of
continuous functions defined on $I^\alp$, indexed by
$\Psi\in\FF_q$,
\[\FF_{C}:=\{c_\Psi(x): \Psi\in\FF_q\}\quad\hbox{and}\quad
\FF_{D}:=\{d_\Psi(x): \Psi\in\FF_q\}.\] We will assume that the
classes $\FF_{C}$ and $\FF_{D}$ are relatively compact with
respect to the sup-norm topology on $I^\alp$, which by the
Arzel\`a-Ascoli theorem is equivalent to these classes being
uniformly equicontinuous and uniformly bounded on
$I^\alp$.\vskip5pt

\noindent We can now state, in Theorem \ref{theo_W} below, the
technical result that will be instrumental in the proof of Theorem
\ref{theo_reg} (see Section \ref{section_final_proof}).
\begin{theo}\label{theo_W}
Let $\{h'_n\}_{n\geq1}$ and $\{h''_n\}_{n\geq1}$ two sequences of
positive constants fulfilling the hypotheses $(H.1$-$2$-$3)$, with
$0<h'_n\leq h''_n<1$. Under the hypotheses $(F.I$-$2$-$IV$-$V)$
and $(K.1$-$2$-$3)$, we have with probability one,
\begin{equation}
\lim_{\ninf} \suph \sup_{\Psi\in\FF_q} \frac{\sup_{{\bx}\in I} \pm
W_{n,h}({\bx},\Psi)}{\sqrt{2n h^d \log (1/h^d)}} =
\sup_{\Psi\in\FF_q} \sigma (\Psi), \label{w}
\end{equation}
\mbox{where}
\begin{equation}
\sigma^{2}(\Psi )=\sup_{{\bx}\in I}
{\E}\Big\{\big(c_\Psi({\bx})\Psi ({\bf Y})+d_\Psi({\bx})\big)^{2}
\,\big|\,
\boldsymbol{X}={\bx}\Big\}f_{\boldsymbol{X}}({\bx})\int_{\R^d}K^{2}({\bf
t})d{\bf t}. \label{ww}
\end{equation}
\end{theo}
\noindent The proof of Theorem \ref{theo_W} is presented in the
following Section \ref{vvvv3}.

\subsection{Proof of Theorem \ref{theo_W}}\label{vvvv3} \noindent
We will mostly borrow the arguments of \cite{EM00}, \cite{EM05}
and \cite{DM04}. First, we will consider the case where
$\mathcal{F}_q$ is reduced to $\{\Psi\}$, for a given real valued,
measurable and uniformly bounded function $\Psi$ defined on $\R^q$
(in such a way that $(F.IV)$ and $(F.V)$ are automatically
fulfilled). Next, in Section \ref{unif_psi}, we will show how to
extend this primary result to a more general class $\FF_q$.

\subsubsection{The case where $\FF_q=\{\Psi\}$} \noindent Assume
that $\FF_q=\{\Psi\}$, where $\Psi$ is a given real valued,
measurable and uniformly bounded function defined on $\R^q$.
Without loss of generality, we will assume that $\kappa = 1$ in
$(K.II)(i)$. Moreover, we will suppose that $\Psi$ is such that
$\sigma(\Psi)>0$ (with $\sigma(\Psi)$ as in $(\ref{ww})$). In this
paragraph, our aim is to establish the following result.
\begin{theo}\label{theo_Wborne}
Let $\{h'_n\}_{n\geq1}$ and $\{h''_n\}_{n\geq1}$ be two sequences
of positive constants fulfilling the hypotheses $(H.1$-$2$-$3)$,
with $0<h'_n\leq h''_n<1$. Under the hypotheses $(F.I$-$2)$,
$(K.1$-$2$-$3)$, we have almost surely
\begin{equation}
\lim_{\ninf} \bigg\{ \suph \Big\{2nh^d\log (1/h^d)\Big\}^{-1/2}
\sup_{{\bx}\in I}\pm W_{n, h} ({\bx},\Psi ) \bigg\}= \sigma
(\Psi)\label{w2}.
\end{equation}
\end{theo}
\noindent The proof of Theorem \ref{theo_Wborne} will be split
into a lower bound part and an upper bound part. \vskip5pt

\noindent {\bf Lower bound part} \quad In this part, we do not
work uniformly over $\inth$. Instead, the result is stated for
$h=h'_n$. It is straightforward that, whenever Proposition
$\ref{borninf}$ below holds, such is also the case when taking the
supremum over $[h'_n,h''_n]$.
\begin{prop}\label{borninf}
Under the assumptions of Theorem \ref{theo_Wborne}, we have with
probability one,
\begin{equation}
\liminf_{n\rightarrow \infty}\Big\{\sup_{{\bx}\in I}\pm W_{n,
h_n'}({\bx},\Psi )/ \sqrt{2n h^{\prime d}_{n}\log (1/h^{\prime
d}_{n})}\Big\}\geq\sigma (\Psi). \label{lit}
\end{equation}
\end{prop}
\demo\ This result is the multidimensional extension of
Proposition 3 in \cite{EM00}. A close look into their proofs
reveals that our "extension" follows from exactly the same lines.
We omit the details of these book-keeping arguments for the sake
of conciseness. \findem \vskip5pt

\noindent {\bf Upper bound part} \quad Here we claim that, under
the assumptions of Theorem \ref{theo_Wborne}, we have, with
probability one,
\begin{equation}
\forall \e >0 \, , \limsup_{\ninf} \Bigg\{\suph
\frac{\sup_{{\bx}\in
    I}|W_{n, h}({\bx},\Psi)|}{\sqrt{2nh^d\log (1/h^d)}}\Bigg\} \leq(1+2\e)\sigma(\Psi).
\label{upper}
\end{equation}
\noindent For any real function  $\varphi$ defined on a set $B$,
we use the notation $\|\varphi\|_B = \sup_{{\bx}\in
B}|\varphi({\bx})|$, and in the particular case where $B=\R^m,
m\geq 1$, we will write $\|\varphi\|_B =\|\varphi\|$. Moreover,
for any class ${\cal G}$ of measurable functions, we will use the
notation
$$
\big\|n^{1/2}\alpha_{n}\big\|_{{\cal G}}:=\sup_{g\in {\cal
G}}\big| n^{1/2}\alpha_{n}(g)\big|.
$$
First note that since $\Psi$ is bounded,  there exists  a constant
$0<M_\Psi<\infty$ such that $\Vert\Psi\Vert\leq M_\Psi$. \vskip3pt

\noindent Now, fix $\e>0$ in $(\ref{upper})$ and introduce some
constants $\gam>0$, $\del\in(0,\alp/4)$, and $\lab\in(1,
(1+2\gamma)^{1/d})$, which will be expressed  in function of
$\e>0$ latter on. For any integer $k\geq 0$, set
\begin{equation}
n_k:= \lfloor(1+\gam)^k\rfloor\label{nkgam1},
\end{equation}
where $\lfloor u\rfloor\leq u<\lfloor u\rfloor+1$ is the integer
part of $u$.\vskip5pt

\noindent For each $k\geq 1$ and $n$ such that $n_{k-1}\leq n\leq
n_k$, consider the interval
\[[h'_{n_k},h''_{n_{k-1}}]
\supset [h'_n,h''_n].\]
Next, set
\begin{equation}
R_k:=
\Big\lfloor\frac{\log(h''_{n_{k-1}}/h'_{n_k})}{\log\lambda}\Big\rfloor+1\label{Rk1},
\end{equation}
and consider the following partitioning of
$[h'_{n_k},h''_{n_{k-1}}]$
\begin{eqnarray}
h'_{n_k,R_k}:= h''_{n_{k-1}}\quad\hbox{and}\quad
h^{\prime}_{n_k,l}:=&\lab^{l} h'_{n_k}\quad\hbox{for}\quad
l=0,\ldots, R_k-1.\label{discret11}
\end{eqnarray}
For each integer $0\leq l\leq R_k$, we now include $I$ in a union
of $J_{l}$ \emph{non-overlapping} hypercubes, with sides of length
$d^{-1/2}{\del h^{\prime}_{n_k,l}}\;$. These hypercubes are
denoted by
\begin{equation}
 \Gam_{k,l,j}:=\Lac x_{k,l,j}+]\,0,d^{-1/2}\del
h^{\prime}_{n_k,l}\;]^d\Rac \quad\hbox{for}\quad 1\le j\le
J_{l},\label{Gamklj}
\end{equation}
and satisfy the following relations, for any integers $k\geq 1$
and $l=0,\ldots, R_k$ (keep in mind that $0<\del<\alp/4$ and
$h'_n<1$ for $n\geq 1$).
\begin{equation}
I\subset
\bigcup_{j=1}^{J_{l}}\Gam_{k,l,j}\subset{I^{\alp/2}}.\label{recouvrement1}
\end{equation}
Since the construction (\ref{Gamklj}) implies that the
$\Gam_{k,l,j}$, $1\le j\le J_{l},$ do not overlap, we can deduce
that there exists a constant $C:=C(\del)$, depending only on
$\del>0$ and  $I$, such that
\begin{equation}
J_{l}\le \frac{C}{h^{\prime d}_{n_k,l}},\quad\hbox{for}\quad k\geq
1\quad\hbox{and}\quad 0\le l\le R_k \label{T2b}.
\end{equation}
In the sequel, we set
\begin{equation}
N_k:=\{n_{k-1}+1,\ldots ,n_k\}\label{Nkgam1},
\end{equation}
if $n_{k-1}< n_k$, and $N_k:=\emptyset$ if $n_{k-1}=n_k$. Observe
that for any initial choice of $\gam>0$, the set $N_k$ is
non-empty provided the integer $k$ is large enough.\vskip5pt

\noindent Our proof begins with the following decomposition, with
 $k\geq 1$ large enough to ensure that $N_k\not=\emptyset$.
\begin{equation}
\begin{split}
\p_k :=\ &\p\lcro \maxn \suph \frac{\sup_{{\bx}\in
I}|W_{n,h}({\bx},\Psi)|}{\sqrt{2nh^d\log(1/h^d)}} \geq
(1+2\e)\sigma(\Psi)
\rcro \\
=\ &\p\lcro \maxn \suph \frac{\sup_{{\bx}\in
I}|\alpha_n(\eta_{{\bx},h,\Psi})|}{\sqrt{2 h^d\log(1/h^d)}} \geq
(1+2\e)\sigma(\Psi)
\rcro \\
\leq\ &\p\lcro \maxn \suphk \frac{\sup_{{\bx}\in
I}|\alpha_n(\eta_{{\bx},h,\Psi})|}{\sqrt{2 h^d\log(1/h^d)}} \geq
(1+2\e)\sigma(\Psi)
\rcro \\
\leq\ &\p\lcro \maxn \maxlrk \maxjjl
\frac{|W_{n,\hkl}(\xklj,\Psi)|}{\sqrt{2n_k \hkld\log(1/\hkld)}}
\geq (1+\e)\sigma(\Psi)
\rcro \\
+ &\p\lcro \maxn \maxoscillolj \suposcillohx \Big|
\frac{\alpha_n(\eta_{{\bx},h,\Psi})}{\sqrt{2h^d\log(1/h^d)}}
\\&\;\;\;-\frac{W_{n,\hkl}(\xklj,\Psi)} {\sqrt{2n_k
\hkld\log(1/\hkld)}}\Big| \geq \e\sigma(\Psi)
\rcro \\
=:\ &\p_{1,k} + \p_{2,k}. \label{decompPk}
\end{split}
\end{equation}
In view of (\ref{decompPk}), our aim is now to prove that
$\p_{1,k}$ and $\p_{2,k}$ are widely summable in $k$ in order to
apply the Borel-Cantelli lemma.\vskip5pt

\noindent {\it Evaluation of $\p_{1,k}$ : partitioning.} Keep in
mind the definition (\ref{hz}) of the functions
$\eta_{\bx,h,\Psi}$. Fix $\gam
>0$ and, for every $k\geq 1$, $0\leq l\leq R_k$, $1\leq j\leq J_l$, and
all $({\bf u},{\bf v})\in \R^d\times\R^q$, define
\begin{equation}
g_{k,l,j}({\bf u},{\bf v}):=\eta_{\xklj,\hkl,\Psi}({\bf u},{\bf
v}). \label{mas2}
\end{equation}
Further introduce, for every $ k\geq 1 $, the class of functions
defined on $\R^d\times\R^q$,
$$
{\cal G}_{k}:=\{g_{k,l,j}:0\leq l\leq R_k, 1\leq j\leq J_l\}.
$$
In view of $(\ref{hz})$ and $(\ref{mas2})$, observe that, for
every $0\leq l\leq R_k$, $1\leq j\leq J_l$ and  all $\inthk$,
\begin{equation}
\Vert g_{k,l,j}\Vert+\Vert \eta_{{\bx},h,\Psi}\Vert \leq 2\big\{
\Vert c_\Psi\Vert\times \Vert\Psi\Vert  +\Vert d_\Psi\Vert
\big\}\Vert K\Vert =:M_1, \label{bd}
\end{equation}
where $\Vert\Psi\Vert\leq M_\Psi<\infty$ by assumption.
%où $\mathcal{{\bf Y}}_J$ est un compact de $\R^q$ tel que ${\bf
%Z}\1\{\boldsymbol{X}\in J\}\in \mathcal{{\bf Y}}_J$ (l'existence
%d'un tel compact est assurée par l'hypothèse $(F.III)$); $\Psi$
%étant supposée bornée sur tout compact de $\R^q$, on a
%$\Vert\Psi\Vert_{\mathcal{{\bf Y}}_J}<\infty$.
\begin{prop}
Assume  the conditions of Theorem $\ref{theo_Wborne}$ are
satisfied. If $(\ref{Gamklj})$ is fulfilled with $0<\delta<\alp/4
$, then we have, almost surely, for all $\e>0$,
\begin{equation}
\limsup_{k \rightarrow \infty}\Bigg\{ \max_{\substack{0\leq l
  \leq R_k\\ 1 \leq j \leq J_l}}
\frac{\maxn|n^{1/2}\alpha_{n}(g_{k,l,j})|}{\sqrt{2 n_k\hkld\log
(1/\hkld)}}\Bigg\}
 \leq (1+\e)\sigma(\Psi). \label{claim1}
\end{equation}
\end{prop}
\noindent\demo\ From $(\ref{T2b})$, we have
\begin{equation}
\begin{split}
\p_{1,k}&\leq \sum_{l=0}^{R_k}\frac{C(\del)}{\hkld}\maxjjl \p
\LCro \frac{\maxn|n^{1/2}\alpha_{n}(g_{k,l,j})|}{\sqrt{2
n_k\hkld\log (1/\hkld)}}\geq (1+\e)\sigmapsi\RCro \\
& =: \sum_{l=0}^{R_k}\frac{C(\del)}{\hkld}\maxjjl \p_{1,k,l,j}.
\end{split} \label{p1klj}
\end{equation}
Arguing as in the proof of Proposition 4.1 in \cite{DM04} (see
especially the proof of their statement (4.11)), it can be shown
that
\begin{equation}
\max_{\substack{0\leq l \leq R_k\\1\leq j \leq J_l}} \maxn \ {\rm
Var}\Big( g_{k,l,j}(\boldsymbol{X},{\bf Y})\Big)\leq \sigma
^{2}(\Psi )(1+\e)\hkld. \label{majvar}
\end{equation}
Thus, recalling the bound (\ref{bd}) on $g_{k,l,j}$, the maximal
version of Bernstein's inequality (see, e.g., Lemma 2.2 in
\cite{EM96}) when applied to the variables
$U_i=g_{k,l,j}(\boldsymbol{X}_i,{\bf Y}_i)-\E
\,g_{k,l,j}(\boldsymbol{X}_i,{\bf Y}_i),\ i=1,\ldots,n$, yields
under $(H.2)$,
\begin{eqnarray}
\p_{1,k,l,j}%&=&{\p}\left\{ \maxn |n^{1/2}\alpha
%  _{n}(g_{k,l,j})|\geq \sigma (\Psi)(1+\e)
%  \sqrt{2n_k \hkld \log (1/\hkld)}\right\} \\[0.2cm]
%&\leq &2\exp \Bigg( -\frac{2\sigma(\Psi )^{2}(1+\e)\hkld\log
%(1/\hkld)n_k}{2n_k\sigma^{2}(\Psi
%)\hkld+ \frac{\textstyle 2M_1\sigma(\Psi )}{\textstyle 3}\sqrt{2n_k\hkld \log (1/\hkld)}}\Bigg), \\[0.2cm]
\leq 2h_{n_k,l}^{\prime\ d(1+\e/2)}. \label{prem_bound}
\end{eqnarray}
Moreover, in view of $(\ref{Rk1})$ we have $R_k\leq
1+[\log(h''_{n_{k-1}}/h'_{n_k})]/\log(\lab)$. Therefore, combining
(\ref{prem_bound}) with $(\ref{p1klj})$, we get
\begin{equation}
\p_{1,k} \leq\ %2C(\del)\sum_{l=0}^{R_k}h_{n_k,l}^{\prime d\e /2}
%%&=\ 2C(\del)h_{n_k}^{\prime d\e /2}\sum_{l=0}^{R_k}\lab^{ld\e /2}\\
%\leq\ \frac{2C(\del)h_{n_k}^{\prime d\e/2}}{\lab^{d\e /2}-1}
%\lab^{d(R_k+1)\e/2}\leq
\frac{2C(\del)\lab^{d\e}}{\lab^{d\e
/2}-1}h_{n_{k-1}}^{\prime\prime\ d\e/2}. \label{finP1k}
\end{equation}
Now, observing that, under $(H.3)$, $\forall \gamma
>0, \ \sum_{k\geq 1}h_{n_{k-1}}^{\prime\prime\ \gamma} < \infty$, the proof of
(\ref{claim1}) is completed by making use of the Borel-Cantelli
lemma. $\sqcup\!\!\!\!\sqcap$ \vskip5pt

\noindent {\it Evaluation de $\p_{2,k}$ : evaluation of the
oscillations. } Set $\e_1 = \e/2$ and $$B_{n,k,l,h} = \bigg|
\sqrt{\frac{2n_k\hkld\log(1/\hkld)}{2nh^d\log(1/h^d)}}-1\bigg|.$$
It holds that
\begin{eqnarray}
\p_{2,k}&=& \p\lcro \maxn \maxoscillolj \suposcillohx \Big|
\frac{W_{n,h}({\bx},\Psi)}{\sqrt{2nh^d\log(1/h^d)}}\nonumber\\
&&
\qquad\qquad\qquad\qquad\qquad\qquad\qquad\quad-\frac{W_{n,\hkl}(\xklj,\Psi)}
{\sqrt{2n_k \hkld\log(1/\hkld)}}\Big| \geq 2\e_1\sigma(\Psi)
\rcro \nonumber\\
&\leq& \p \lcro \maxn \maxoscillolj \!\!\suposcillohx \!\!\!\!\!
\bigg| \frac{W_{n,h}({\bx},\Psi)-W_{n,\hkl}(\xklj,
\Psi)}{\sqrt{2n_k\hkld\log(1/\hkld)}} \bigg|\geq \e_1\sigmapsi \rcro \nonumber\\
&&\!\!\!\!\!+\ \p \lcro \maxn \maxoscillolj \suposcillohx \!\!\!
B_{n,k,l,h}\bigg|\frac{W_{n,h}({\bx},\Psi)}{\sqrt{2n_k\hkld\log(1/\hkld)}}\bigg|\geq
\e_1\sigmapsi\rcro \nonumber\\
&=:& \p_{2,1,k}+\p_{2,2,k}. \label{decompP2k}
\end{eqnarray}
First consider $\p_{2,2,k}$, and note that, for every $n\in N_k$
and all $\inthkl$,
\begin{equation}
B_{n,k,l,h} \leq \sqrt {\frac{\log(1/\hkld)}
{\log(1/\lab^d\hkld)}}-1\bigg| \sqrt{\frac{\lab^d n_k}{n}}+ \bigg|
\sqrt {\frac{\lab^d n_k}{n}} - 1 \bigg|.
\end{equation}
Some algebra enable to state
\begin{equation}
\liminf_{k\tend\infty} \min_{n\in N_k} \inf_{\inth} \e_1\sigmapsi
B_{n,k,l,h}^{-1}\geq 2\Lpar1+\sqrt{\frac{2}{A_2}}\Rpar
D_1(\nu)\sigmapsi, \label{liminfepsiBklh}
\end{equation}
where $A_2$ and $D_1(\nu)$ are the constants involved in Fact
$\ref{factvarron_reg}$ (see the Appendix).\\
Now set, for $0\leq l\leq R_k-1$ and $1\leq j\leq J_l$,
\[\p_{2,2,k,l,j}:= \p \lcro \maxn {\displaystyle\suposcillohx
B_{n,k,l,h}\bigg|\frac{W_{n,h}({\bx},\Psi)}{\sqrt{2n_k\hkld\log(1/\hkld)}}}\bigg|\geq
\e_1\sigmapsi\rcro.\] In view of (\ref{decompP2k}), we have
\begin{eqnarray}
\p_{2,2,k} \leq \sum_{l=0}^{R_k-1}\sum_{j=1}^{J_l} \p_{2,2,k,l,j}.
\label{def_p22klj}
\end{eqnarray}
Keep in mind the definition (\ref{Gkprime}) of $\GG'$ and
introduce the class of functions
\begin{equation}
\overline{\FF}_{k,l,j} := \lac (c_\Psi({\bx})+d_\Psi({\bx})\Psi)K
\lpar \frac{\cdot - {\bx}}{h}\rpar, {\bx}\in \Gamklj, \inthkl \rac
\subset \GG'.  \label{Fbarklj}
\end{equation}
It is easy to check that $\overline{\FF}_{k,l,j} $ is a pointwise
measurable class of bounded functions admitting a polynomial
uniform covering number, for all $0\leq l\leq R_k$ and $1\leq
j\leq J_l$. Moreover, arguing as above (see $(\ref{majvar})$), it
can be shown that uniformly over $\inthkl$ and ${\bx}\in \Gamklj$,
\begin{equation}
\begin{split}
{\rm Var}\lcro (c_\Psi({\bx})+d_\Psi({\bx}))\Psi({\bf
Y}))K\lpar\frac{\boldsymbol{X}-{\bx}}{h}\rpar
\rcro &\leq \lab^d \sigmapsi^2\hkld + o(\hkld)\\
%&\leq 2 \sigmapsi^2\hkld + o(\hkld)\\
&\leq 4 \sigmapsi^2\hkld,\label{var1}
\end{split}
\end{equation}
for $k$ large enough. Further observe that $\|g\| \leq M_1$ for
all $g\in\overline{\FF}_{k,l,j}$, with $M_1$ as in (\ref{bd}).
Therefore, we can apply Fact \ref{factvarron_reg} with $\tau =
2\sigmapsi$ and $\rho=\tau\sqrt{2/A_2}$. When combined with
$(\ref{liminfepsiBklh})$, this yields
\begin{equation}
\begin{split}
\p_{2,2,k,l,j}&\leq \p \lcro \maxn \suposcillohx
\bigg|\frac{W_{n,h}({\bx},\Psi)}
{\sqrt{2n_k\hkld\log(1/\hkld)}}\bigg|
\geq 2\Big(1+\sqrt{\frac{2}{A_2}}\Big)D_1(\nu)\sigmapsi\rcro \\
&\leq \p\lcro \max_{1\leq n\leq n_k}
\|n^{1/2}\alpha_n\|_{\overline{\FF}_{k,l,j}}\geq
D_1(\nu)(\tau+\rho)
\sqrt{2n_k\hkld\log(1/\hkld)} \rcro \\
%&\leq 4\exp \Lcro -A_2\big(\frac{\rho}{\tau}\big)^2\log(1/\hkld)\Rcro\\
&\leq 4h_{n_k,l}^{\prime\ 2d}.\label{bien_utile_finalement}
\end{split}
\end{equation}
In view of (\ref{def_p22klj}) and (\ref{bien_utile_finalement}),
and arguing as in $(\ref{finP1k})$, it follows that
\begin{equation}
\p_{2,2,k}\leq \frac{4C(\del)}{\lab^d-1}h^{\prime \prime\
d}_{n_{k-1}}. \label{finP22k}
\end{equation}
\vskip8pt

\noindent Turning our attention to $\p_{2,1,k}$, set, for $0\leq
l\leq R_k-1$ and $1\leq j\leq J_l$,
\[\p_{2,1,k,l,j}:= \p \lcro \maxn \displaystyle \suposcillohx
\bigg| \displaystyle \frac{W_{n,h}({\bx},\Psi)-W_{n,\hkl}(\xklj,
\Psi)}{\sqrt{2n_k\hkld\log(1/\hkld)}} \bigg|\geq \e_1\sigmapsi
\rcro.\] In view of (\ref{decompP2k}), it follows that
\begin{eqnarray}
\p_{2,1,k} \leq \sum_{l=0}^{R_k-1}\sum_{j=1}^{J_l} \p_{2,1,k,l,j}.
\label{def_p21klj}
\end{eqnarray}
Next consider the classes of functions
\begin{equation}
\overline{\FF}'_{k,l,j} := \lac g_{k,l,j}-\eta_{{\bx},h,\Psi},
{\bx}\in \Gamklj, \inthkl \rac\subset\GG', \label{Fbarprimeklj}
\end{equation}
with $g_{k,l,j}$ defined as in $(\ref{mas2})$. To evaluate
$\p_{2,1,k,l,j}$, our aim is to apply Fact \ref{factvarron_reg}
once again. First note that $\|g\|\leq M_1$ for all $g\in
\overline{\FF}'_{k,l,j}$ (with $M_1$ still as in $(\ref{bd})$).
Now we shall get a suitable upper-bound for $\sup_{g\in
\overline{\FF}'_{k,l,j}}{\rm Var}(g(\boldsymbol{X},{\bf Y}))$. Set
$\beta:=\Vert \Psi \Vert^{2}+1<\infty $ and introduce, for any
real valued function $\phi$ defined on $I^\alp$ and all
$\delta\geq0$,
\begin{equation}
\omega_{\phi }(\delta ):=\sup \Big\{|\phi ({\bx})-\phi
(\boldsymbol{y})|:\|{\bx}-\boldsymbol{y}\|\leq \del ,\hbox{ and
}{\bx},\boldsymbol{y}\in I^\alp\Big\}. \label{modcont}
\end{equation}
Making use of the classical inequality $(a+b)^{2}\leq
2(a^{2}+b^{2})$, it can be derived, from $(\ref{hz})$,
$(\ref{Gamklj})$ and $(\ref{mas2})$, that, for $k\geq 1$,
$\intlrk$ and $\intjjl$,
\begin{equation}
\begin{split}
&{\E}\Big\{\Big(g_{k,l,j}(\boldsymbol{X},{\bf
Y})-\eta_{n,{\bx},h}(\boldsymbol{X},{\bf Y})\Big)^2\Big\} =
{\E}\Big\{\Big(\eta_{n_k,\xklj,\hkl}(\boldsymbol{X},{\bf Y})-
\eta_{n,{\bx},h}(\boldsymbol{X},{\bf Y})\Big)^2\Big\}\\
%=&\ {\E} \Big\{\Big(\big(c_\Psi(\xklj)\Psi({\bf Y})+
%d_\Psi(\xklj)\big)K\Big({\xklj-\boldsymbol{X}\over
%\hkl}\Big)-(c_\Psi({\bx})
%\Psi({\bf Y})+d_\Psi({\bx})\big)K\Big({{\bx}-\boldsymbol{X}\over h}\Big)\Big)^2\Big\}\\
\leq &\
2{\E}\Big\{\Big(\big(c_\Psi({\bx}_{k,l;j})-c_\Psi({\bx})\big)
\Psi({\bf
Y})+d_\Psi(\xklj)-d_\Psi({\bx})\Big)^2K\Big(\frac{\xklj-\boldsymbol{X}}
{\hkl}\Big)^2 \Big\}\\
+&\ 2{\E}\Big\{ \big(c_\Psi({\bx})\Psi({\bf
Y})+d_\Psi({\bx})\big)^{2}\Big(
K\Big(\frac{\xklj-\boldsymbol{X}}{\hkl}\Big) -
K\Big(\frac{{\bx}-\boldsymbol{X}}
{h}\Big) \Big)^{2}\Big\} \\
\leq\ &4\beta\Big\{ \omega_{c_\Psi}^{2}(\delta \hkl)\vee
\omega_{d_\Psi}^{2}(\delta \hkl)\Big\}  \,{\E}\Big\{ K \Big(
\frac{\xklj-\boldsymbol{X}}{\hkl}\Big) ^{2}\Big\} \\
+&\ 4\beta\Big\{ \Vert c_\Psi\Vert _{I^\alp}^{2} \vee \Vert
d_\Psi\Vert _{I^\alp}^{2}\Big\}{\E}\Big\{ \Big( K\Big(
\frac{\xklj-\boldsymbol{X}}{\hkl}\Big)
-K\Big(\frac{{\bx}-\boldsymbol{X}}{h}\Big) \Big) ^{2}\Big\}\\
=:&\ (I)+(I\kern-0.18em I). \label{decompVaroscillo}
\end{split}
\end{equation}
Set
\begin{equation}
\widetilde{\Gam}_{k,l,j}:=\Lac \xklj+[-\hkl/2, \hkl/2]^d\Rac.
\label{Gamtildklj}
\end{equation}
Assuming as before that $\kappa=1$ in $(K.II)(i)$, we obtain
\begin{eqnarray*}
K\Big(\frac{\xklj-\boldsymbol{X}}{\hkl}\Big) ^{2}\leq  \Vert
K\Vert ^{2} {\1}{\big\{ \boldsymbol{X}\in \widetilde{\Gam}_{k,l,j}
\big\}}.
\end{eqnarray*}
So we have, for $k$ large enough, and for every $\intlrkm$ and
$\intjjl$,
\begin{eqnarray}
(I)&\leq&4\beta \Big\{\omega_{c_\Psi}^{2}(\delta \hkl)\vee \omega
_{d_\Psi}^{2}(\delta \hkl)\Big\}\,\|K\|^{2}\,\p(\boldsymbol{X}\in\widetilde{\Gam}_{k,l,j}) \nonumber \\[0.2cm]
&\leq&4\beta \Big\{\omega_{c_\Psi}^{2}(\delta \hkl)\vee \omega
_{d_\Psi}^{2}(\delta \hkl)\Big\}
\,\|K\|^{2}\,\|f_{\boldsymbol{X}}\|_{I^\alp}\hkld. \label{I1}
\end{eqnarray}
In order to bound $(I\kern-0.18em I)$, we will employ the
assumption $(K.I)$. Let $B_1$ and $B_2$ be the two functions
respectively defined on $\R$ and $[1,+\infty)$ by
\begin{eqnarray}
B_1(\delta)&=&\sup_{|{\bf
u}|\leq\delta}\int_{\R^d}\big(K(\bx)-K(\bx + {\bf
u})\big)^2d\bx, \label{def_A1}\\
B_2(\gamma)&=&\sup_{1/\gamma\leq\lambda\leq\gamma}\int_{\R^d}\big(K(\lab\bx)-K(\bx)\big)^2d\bx.
\label{def_A2}
\end{eqnarray}
The assumption $(K.I)$ ensures that
\begin{eqnarray}
\lim_{\delta\tend0}B_1(\delta)=0 \quad\mbox{  and  }\quad
\lim_{\lab\tend 1}B_2(\lab)=0. \label{lim_A1A2}
\end{eqnarray}
Let $\Delta:=\int_{\R^d}(K((\bx-{\bf t})/h)-K((\xklj-{\bf
t})/\hkl))^2d{\bf t}$. Setting  ${\bf u}=(\xklj-{\bf t})/\hkl$, we
obtain, from the definitions (\ref{discret11}), (\ref{Gamklj}),
(\ref{def_A1}) and (\ref{def_A2}), and the Cauchy-Schwarz
inequality,
\begin{equation} \nonumber
\begin{split}
\hkld\Delta &=\int_{\R^d}\Big[K\Big(\frac{\bx-{\xklj}}{h}+{\bf
u}\frac{\hkl}{h}\Big)-K({\bf u})\Big]^2d{\bf u}\\
%&=\int_{\R^d}\Big[K\Big(\frac{\bx-{\xklj}}{h}+{\bf
%u}\frac{\hkl}{h}\Big)-K\Big({\bf u}\frac{\hkl}{h}\Big)+K\Big({\bf
%u}\frac{\hkl}{h}\Big)-K({\bf u})\Big]^2d{\bf u}\\
&=\int_{\R^d}\Big[K\Big(\frac{\bx-{\xklj}}{h}+{\bf
u}\frac{\hkl}{h}\Big)-K\Big({\bf u}\frac{\hkl}{h}\Big)\Big]^2d{\bf
u}+\!\int_{\R^d}\! \Big[K\Big({\bf u}\frac{\hkl}{h}\Big)-K({\bf
u})\Big]^2\\
&\quad + 2\int_{\R^d} \Big[K\Big(\frac{\bx-{\xklj}}{h}+{\bf
u}\frac{\hkl}{h}\Big)-K\Big({\bf
u}\frac{\hkl}{h}\Big)\Big]\Big[K\Big({\bf
u}\frac{\hkl}{h}\Big)-K({\bf
u})\Big]d{\bf u}\\
%&=: \Delta_1+\Delta_2+\Delta_3.
&\leq 2B_1\Big(\frac{2\delta}{\sqrt{d}}\Big) + B_2(\lab) +
2\sqrt{2B_1\Big(\frac{2\delta}{\sqrt{d}}\Big) B_2(\lab)}.
\end{split}
\end{equation}
Therefore, we have $\Delta\leq
\Big(\sqrt{2B_1\big(\frac{2\delta}{\sqrt{d}}\big)} +
\sqrt{B_2\big({\lab}\big)} \Big)^2$, in such a way that
\begin{eqnarray*}
(I\kern-0.18em I)\leq 4\Big\{ \Vert c_\Psi\Vert _{I^\alp}^{2} \vee
\Vert d_\Psi\Vert _{I^\alp}^{2}\Big\} \beta
\|f_{\boldsymbol{X}}\|_{I^\alp}
\Big(\sqrt{2B_1\Big(2\delta/\sqrt{d}\Big)} +
\sqrt{B_2\big({\lab}\big)} \Big)^2 \hkld.
\end{eqnarray*}
From $(\ref{decompVaroscillo})$ and $(\ref{I1})$, setting
$B=4\|f_{\boldsymbol{X}}\|_{I^\alp}\beta\Lac\|K\|^2 + \big(
\|c_\Psi\|_{I^\alp}^2\vee\|d_\Psi\|_{I^\alp}^2\big)|K|_v^2\Rac$,
it holds that
\begin{eqnarray*}
{\rm Var}(g_{k,l,j}(\boldsymbol{X},{\bf
Y})-\eta_{{\bx},h,\Psi}(\boldsymbol{X},{\bf Y})) &\leq& B\Big\{(
\omega_{c}^{2}(\del\hkl)\vee \omega_{d}^{2}(\del\hkl)\\
&&\quad \quad+\Big(\sqrt{2B_1\Big(2\delta/\sqrt{d}\Big)} +
\sqrt{B_2\big({\lab}\big)} \Big)^2 \Big\}\hkld.
\end{eqnarray*}
Now set for $k\geq1$, $\intlrkm$ and $\intjjl$,
$$
\sigma _{k,l,j}^{2}(\Psi)=\sup \Big\{{\rm
Var}(g(\boldsymbol{X},{\bf Y})):g \in
\overline{\FF}'_{k,l,j}\Big\}.
$$
By selecting $\del>0$ small enough and $\lab>1$ close enough to 1,
the continuity of the functions $c_\Psi$ and $d_\Psi$ when
combined with (\ref{lim_A1A2}) implies that, for $k$ large enough,
\begin{equation}
\sigma _{k,l,j}^{2}(\Psi) \leq \frac{\e_1^2\sigmapsi^2}{[(1+
\sqrt{2/A_2})D_1(\nu)]^2}\hkld,
\end{equation}
where $A_2$ and $D_1(\nu)$ are the constants involved in Fact
$\ref{factvarron_reg}$. Moreover, arguing as above, it is easy to
check that $\overline{\FF}'_{k,l,j}\subset\GG'$ is a pointwise
measurable class of bounded functions with a polynomial uniform
covering number for every $\intlrkm$ and $\intjjl$. Therefore, we
can apply Fact $\ref{factvarron_reg}$ with
$\tau=\e_1\sigmapsi/[(1+ \sqrt{2/A_2})D_1(\nu)]$ and $\rho=\tau
\sqrt{2/A_2}$, which yields
\begin{equation}
\begin{split}
\p_{2,1,k,l,j} &\leq \p\lcro \max_{1\leq n\leq n_k}
\|n^{1/2}\alpha_n\|_{\overline{\FF}'_{k,l,j}}\geq
D_1(\nu)(\tau+\rho)
\sqrt{2n_k\hkld\log(1/\hkld)} \rcro \\
%&\leq 4\exp \Lcro -A_2\big(\frac{\rho}{\tau}\big)^2\log(1/\hkld)\Rcro\\
&\leq 4h_{n_k,l}^{\prime\ 2d}.\label{bien_utile_finalement2}
\end{split}
\end{equation}
Arguing as in  $(\ref{finP1k})$, it follows from
(\ref{def_p21klj}) and (\ref{bien_utile_finalement2}) that
\begin{equation}
\begin{split}
\p_{2,1,k}%&\leq\ 4C(\del)\sum_{l=0}^{R_k-1}\hkld \\
%&\leq\ 4C(\del)h^{\prime d}_{n_k}\sum_{l=0}^{R_k-1}\lab^{dl} \\
\leq \frac{4C(\del)}{\lab^d-1}h^{\prime \prime\ d}_{n_{k-1}}.
\label{finP21k}
\end{split}
\end{equation}
\vskip5pt

\noindent By combining (\ref{finP22k}) and (\ref{finP21k}) we
conclude under $(H.3)$ that $\sum_{k\geq1} \p_{2,k} < \infty$,
which achieves the evaluation of the oscillations part. The
statement $\sum_{k\geq1} \p_{k} < \infty$ now directly follows
from $(\ref{finP1k})$, completing the proof of (\ref{upper}).
\vskip5pt

\noindent By combining (\ref{upper}) with Proposition
\ref{borninf}, the statement of Theorem \ref{theo_Wborne} is
straightforward.

\subsubsection{The general case}\label{unif_psi}
\noindent To extend the result of Theorem \ref{theo_Wborne} to the
general case, only the "upper bound" part has to be extended. From
the preceding paragraph, it is straightforward that for any finite
subclass ${\mathcal G}_q\subset {\cal F}_q$, we have, with
probability one,
\begin{equation}
\limsup_{n\tend\infty}\suph\frac{\sup_{\Psi\in{\mathcal
G}_q}{\supx}|W_{n,h}({\bx},\Psi)|}{\sqrt{2nh^d\log(1/h^d)}}\leq
\sup_{\Psi\in{\mathcal G}_q}\sigma(\Psi)
\label{Borne_sup_ssclasse_finie}.
\end{equation}
Here we shall show that (\ref{Borne_sup_ssclasse_finie}) can be
extended to the entire class ${\cal F}_q$. The following couple of
Lemmas are directed towards this aim. Set
\[C_{\FF_q}:=\sup\{\|c_\Psi\|_{I^{\alp}}:\Psi\in\FF_q\}
\quad\mbox{ and }\quad
D_{\FF_q}:=\sup\{\|d_\Psi\|_{I^{\alp}}:\Psi\in\FF_q\}.\] Keep in
mind that, by $(F.IV)$, the class $\FF_q$ has a uniformly bounded
envelop function $\Upsilon$, with $\Upsilon({\bf
y})\geq\sup_{\Psi\in\FF_q}\Psi({\bf y})$, ${\bf y}\in\R^q.$
\begin{lem}\label{lem_unif_classe_1}
For all $\e>0$, there exists a finite subclass  ${\mathcal
G}_{q,\e}\subset{\mathcal F}_q$ such that, for all
$\Psi\in{\mathcal F}_q$, for $n$ large enough,
\begin{equation}
\min_{\phi\in{\mathcal G}_{q,\e}}\suphx\frac{1}{h^d}
\E\Big[\big\{c_\Psi({\bx})\Psi({\bf {\bf
Y}})+d_\Psi({\bx})-c_\phi({\bx})\phi({\bf {\bf
Y}})-d_\phi({\bx})\big\}^2 K^2\Big(\frac{{\bf
X}-{\bx}}{h}\Big)\Big] \leq \e. \label{res_lem_unif_classe_1}
\end{equation}
\end{lem}
\demo\ Let $J$ be a compact of $\R^q$. For $\Psi, \phi\in{\mathcal
F}_q$, ${\bx}\in I$ and $\inth$, we have
\begin{eqnarray}
&&\frac{1}{h^d}  \E\Big[\big\{c_\Psi({\bx})\Psi({\bf {\bf
Y}})+d_\Psi({\bx})-c_\phi({\bx})\phi({\bf {\bf
Y}})-d_\phi({\bx})\big\}^2 K^2\Big(\frac{{\bf X} -
{\bx}}{h}\Big)\Big]\nonumber\\
&&\;\; \leq \;\;\frac{1}{h^d} \E\Big[\big\{c_\Psi({\bx})\Psi({\bf
{\bf Y}})+d_\Psi({\bx})-c_\phi({\bx})\phi({\bf {\bf
Y}})-d_\phi({\bx})\big\}^2 \1_{\{|{\bf X}-{\bx}|\leq
h/2\}}\Big]\|K^2\|_\infty\nonumber\\
&&\;\;\leq \;\;\frac{4}{h^d}\E[(C_{\FF_q}\Upsilon({\bf
Y})+D_{\FF_q})^2\1_{\{{\bf Y}\notin
J,{|\bf X}-{\bx}|\leq h/2\}}]\|K^2\|_\infty \nonumber\\
&&\;\;\;\;+\;\; \eta(J) \int_J \{c_\Psi({\bx})\Psi({\bf {\bf
Y}})+d_\Psi({\bx})-c_\phi({\bx})\phi({\bf {\bf Y}})-d_\phi({\bx})\big\}^2d{\bf y},\nonumber\\
&&\;\;=:\;\;(I)+ (I\kern-0.22em I).\nonumber
\end{eqnarray}
where $\eta(J):= \sup_{({\bx},{\bf y})\in I^\alp\times J}f_{{\bf
X},{\bf {\bf Y}}}({\bx},{\bf y}) \|K^2\|_\infty$.\\
To evaluate the term $(I)$, first observe that the Hölder
inequality gives us, for all $s> 1$,
\begin{eqnarray*}
(I)\leq 4\|K^2\|_\infty\tilde{\alp}^{2/s}\|f_{\bf
X}\|_{I^\alpha}^2\supx\P({\bf Y}\notin J|{\bf X}={\bx})^{1-2/s},
\end{eqnarray*}
where $\tilde{\alp}:=\supx \E[(C_{\FF_q}\Upsilon({\bf
Y})+D_{\FF_q})^s|{\bf X}={\bx}]<\infty$. Under our continuity
assumptions, and from Scheffé's lemma, the function
${\bx}\rightarrow\P({\bf Y}\in\cdot|{\bf X}={\bx})$ is continuous
from $I$ to the space of all probability measures on Borelian sets
of $\R^q$ endowed with the topology of weak convergence. Thus, the
set $\{\P({\bf Y}\in\cdot|{\bf X}={\bx}), {\bx}\in I\}$ is compact
in this space, which implies, in view of Prohorov's theorem, that
it is uniformly tight. Consequently, for all $\e>0$, there exists
a compact $J_{\e}\subset\R^q$ such that  $\P({\bf Y}\in
J_{\e}|{\bf X}={\bx})>1-\e/2$. Finally, we obtain that there
exists a compact $J=J_{\e}$ such that, uniformly in ${\bx}\in I$,
\begin{eqnarray*}
(I)\leq\frac{\e}{2}.
\end{eqnarray*}
To evaluate the term $(I\kern-0.22em I)$, we will use the fact
that ${\mathcal F}_q$ is a $VC$ subgraph class of functions, which
ensure that it is totally bounded with respect to $d_{Q_\e}$,
where $Q_\e$ is the uniform law over $J^\e$. We can deduce that
for all $\delta>0$, there exists a finite subclass ${\mathcal
G}_{q,1}\subset{\mathcal F}$ for which
\begin{eqnarray*}
\sup_{\Psi\in{\mathcal F}_q}\min_{\phi\in{\mathcal G}_{q,1}}
\int_{J^\e}(\Psi({\bf y})-\phi({\bf y}))^2d{\bf y} <\delta.
\end{eqnarray*}
Moreover, the relative compactness of the classes $\FF_{{\mathcal
C}}$ and $\FF_{{\mathcal D}}$ implies the existence of finite
subclasses  ${\cal G}_{q,2},{\cal G}_{q,3}\subset{\cal F}_q$ such
that
\begin{eqnarray*}
\sup_{\Psi\in{\mathcal F}_q}\min_{\phi\in{\mathcal
G}_{q,2}}\|c_{\Psi}-c_\phi\|_I \vee \sup_{\Psi\in{\mathcal
F}}\min_{\phi\in{\mathcal G}_{q,3}}\|d_{\Psi}-d_\phi\|_I <\delta.
\end{eqnarray*}
Combining the facts that
\begin{eqnarray*}
\sup_{\Psi\in{\mathcal F}_q}\|c_\Psi\|_{I^\alpha}<\infty
\quad\mbox{ and }\quad\sup_{\Psi\in{\mathcal
F}_q}\int_{J^\e}\Psi({\bf y})^2d{\bf y}<\infty,
\end{eqnarray*}
and choosing $\delta>0$ small enough, we have
\begin{equation*}
\sup_{\Psi\in{\mathcal F}_q} \min_{\phi_1, \phi_2, \phi_3} \supx
\int_{J^\e}\{c_\Psi({\bx})\Psi({\bf {\bf
Y}})+d_\Psi({\bx})-c_{\phi_2}({\bx})\phi_1({\bf {\bf
Y}})-d_{\phi_3}({\bx})\big\}^2d{\bf y} \leq \e/(4\eta(J_{\e})),
\end{equation*}
where the minimum is taken over ${\mathcal G}_{q,1}\times{\cal
G}_{q,2}\times{\cal G}_{q,3}$. For any triplet $(\phi_1, \phi_2,
\phi_3)\in{\mathcal G}_{q,1}\times{\cal G}_{q,2}\times{\cal
G}_{q,3}$ for which there exists $\phi\in{\cal F}_q$ such that
\begin{equation*}
\supx \int_{J^\e}\{c_\phi({\bx})\phi({\bf {\bf
Y}})+d_\phi({\bx})-c_{\phi_2}({\bx})\phi_1({\bf {\bf
Y}})-d_{\phi_3}({\bx})\big\}^2d{\bf y} \leq \e/(4\eta(J_{\e})),
\end{equation*}
we select one of these $\phi\in\FF_q$ to construct the subclass
${\mathcal G}_\e$. Applying the triangle inequality, we obtain
\begin{equation*}
(I\kern-0.22em I)\leq\frac{\e}{2},
\end{equation*}
which completes the proof of the lemma. \findem \vskip5pt\noindent
Set $\e>0$ and $n_0$ such that (\ref{res_lem_unif_classe_1}) holds
for all $n\geq n_0$. For any $\Psi,\phi\in\FF_q$, we define,
\begin{eqnarray*}
d^2(\Psi,\phi)&:=& \sup_{n\geq n_0}\suph h^{-d}\supx\E\Big[
\{c_\Psi({\bx})\Psi({\bf Y})+d_\Psi({\bx})\\
&&\;\;\;\;-c_\phi({\bx})\phi({\bf
Y})-d_\phi({\bx})\}^2K^2\Big(\frac{{\bf X}-{\bx}}{h}\Big)\Big].
\end{eqnarray*}

\begin{lem}\label{lem_unif_classe_2}
Under the assumptions of Theorem \ref{theo_W}, there exists an
absolute constant $A$ such that, for all $\e>0$, we have almost
surely,
\begin{equation}
\limsup_{n\tend\infty}\suph\frac{\sup_{d^2(\Psi,\phi)\leq\e}\supx|W_{n,h}(\bx,\Psi)-W_{n,h}(\bx,\phi)|}
{\sqrt{nh^d\log(1/h^d)}}\leq A D_1(\nu) \sqrt{\e}.
\end{equation}
\end{lem}
\demo\  To establish Lemma \ref{lem_unif_classe_2}, we intend to
apply Fact \ref{factvarron_reg}. Consider the classes of functions
\[\tilde\FF(\e,h):=\{\eta_{\bx, h,\Psi}-\eta_{\bx,
h,\phi}:d^2(\Psi,\phi)\leq\e, \bx\in I\},\mbox{ for} \inth,\]
where the functions  $\eta_{\bx, h,\Psi}$ are defined as in
(\ref{hz}). For any integer  $k\geq 1$, set $n_k=2^k$. Further
setting
\begin{eqnarray*}
Q_k:=\max_{n_{k-1}\leq n\leq n_k} \sup_{\inth}
\frac{\sup_{d^2(\Psi,\phi)\leq\e} \sup_{\bx\in
I}|W_{n,h}(\bx,\Psi)-W_{n,h}(\bx,\phi)|}{\sqrt{2nh^d\log(1/h^d)}},
\end{eqnarray*}
we have
\begin{eqnarray*}
Q_k=\max_{n_{k-1}\leq n\leq n_k} \sup_{\inth}
\frac{\|n^{1/2}\alpha_n\|_{\tilde\FF(\e,h)}}{\sqrt{2nh^d\log(1/h^d)}}.
\end{eqnarray*}
Now, in the same way as above, consider the following partitioning
of the interval $[h'_{n_k},h''_{n_{k-1}}]$ $\supset$
$[h'_n,h''_n]$, for $n_{k-1}\leq n\leq n_k$,
\begin{align}
\nono h'_{n_k,R_k}:=& h''_{n_{k-1}}\\
\hkl:=&\lab^l h'_{n_k},\;l=0\ldots R_k-1,\label{discret111}
\end{align}
where $R_k$ satisfy the condition
\begin{equation}
R_k=
\Big\lfloor\frac{\log(h''_{n_{k-1}}/h'_{n_k})}{\log(\lab)}\Big\rfloor+1\label{Rk111}.
\end{equation}
For any $\intlrk$, introduce the classes of functions
\begin{equation}
\tilde{\FF}_{k,l}(\e):=\{\eta_{\bx, h,\Psi}-\eta_{\bx,
h,\phi}:d^2(\Psi,\phi)\leq\e, \bx\in I, h\in[h'_{n_k,l},
h'_{n_k,l+1}\}. \label{def_tildeFF}
\end{equation}
In view of this definition of $\tilde{\FF}_{k,l}(\e)$, we have,
for any function $\upsilon\in \tilde{\FF}_{k,l}(\e)$ and for $k$
large enough,
\begin{equation}
h^{-d}{\rm Var}(\upsilon({\bf X},{\bf Y})) \leq
d^2(\Psi,\phi).\label{majvarunif}
\end{equation}
Thus, from  (\ref{def_tildeFF}) and (\ref{majvarunif}), setting
$\tilde{\sigma}^2_{k,l,\e}:=\sup\{{\rm Var}(\upsilon({\bf X},{\bf
Y})) : \upsilon \in \tilde{\FF}_{k,l}(\e)\},$ we have
$\tilde{\sigma}^2_{k,l,\e}\leq \lambda^{dl}\e h^{'\
d}_{n_k}.$\vskip5pt

\noindent Furthermore, for any function $\upsilon\in
\tilde{\FF}_{k,l}(\e)$, we have $\|\upsilon\|\leq M_1$ (where
$M_1$ is defined as in (\ref{bd})) and we can show that each
subclass $\tilde{\FF}_{k,l}(\e)$ is pointwise measurable and such
that $\tilde{\FF}_{k,l}(\e)\subset \tilde\GG$, where $\tilde\GG$
is a pointwise measurable class of functions admitting a bounded
envelope function and a polynomial uniform covering number (to do
so, we use the same arguments as those used to establish that
$\GG'$ had this property along with the fact that $\FF_q$ is $VC$;
we refer to the proof of the Lemma 5 in \cite{EM00} for more
details). So, the assumptions of Fact \ref{factvarron_reg} are
fulfilled, with $\tau=\e\sqrt{\lambda^{dl}}$ and
$\rho=\tau\sqrt{1/A_2}$. By applying this result, we have, with
$N_k=\{n_{k-1}+1,...,n_k\}$ as usual, and applying the same
techniques as those used to establish
(\ref{bien_utile_finalement}),
\begin{eqnarray*}
\P_{\e,k,l}&:=& \P \Big[\max_{1\leq n \leq
n_k}\frac{n^{1/2}\|\alpha_n\|_{\tilde{\FF}_{k,l}(\e)}}{\sqrt{2n_kh^{\prime
d}_{n_k,l}\log(1/h^{\prime d}_{n_k,l})}}
\geq \sqrt{\lambda^{dl}}(1+\sqrt{1/A_2})D_1(\nu)\e\Big]\\
&\leq& 4 h^{\prime d}_{n_k,l}.
\end{eqnarray*}
Thus,
\begin{eqnarray*}
\P_{\e,k}&:=& \P \Big[\max_{0\leq l\leq R_k-1}\max_{1\leq n \leq
n_k}\frac{n^{1/2}\|\alpha_n\|_{\tilde{\FF}_{k,l}(\e)}}{\sqrt{2n_kh^{\prime
d}_{n_k,l}\log(1/h^{\prime d}_{n_k,l})}}
\geq \sqrt{\lambda^{dl}}(1+\sqrt{1/A_2})D_1(\nu)\e\Big]\\
&\leq& 4 \sum_{l=0}^{R_k-1} h^{\prime d}_{n_k,l}\\
&\leq& \frac{4}{\lambda^d-1}h^{\prime\prime d}_{n_{k-1}},
\end{eqnarray*}
which implies the lemma \ref{lem_unif_classe_2}, \emph{via}
Borel-Cantelli's lemma. \findem
 \vskip5pt\noindent Combining the lemmas \ref{lem_unif_classe_1} and \ref{lem_unif_classe_2} with
(\ref{Borne_sup_ssclasse_finie}), we obtain the uniformity over
the class $\FF_q$, which ends the proof of theorem \ref{theo_W} in
the general case.

\subsection{Proof of Theorem \ref{theo_reg}}\label{section_final_proof}
In this part, we show how the results obtained for the process
$W_{n, h}({\bx}, \Psi)$ can be transposed to nonparametric
functional estimators of the regression, that is how Theorem
\ref{theo_reg} directly follows from Theorem \ref{theo_W}. Towards
this aim, first introduce the following quantities
\begin{eqnarray*}
r_{\psi}(\boldsymbol{x})&=&\int_{\R^q}\psi({\bf
y})f_{\boldsymbol{X},{\bf Y}}(\boldsymbol{x},{\bf y})d{\bf y};\\
f_{\boldsymbol{X},n,h}(\boldsymbol{x})
&=&{\frac{1}{nh^d}}\sum_{i=1}^{n}K\Big(
{\frac{\boldsymbol{x}-\boldsymbol{X}_{i}}{h}}
\Big);  \\
r_{\psi,n,h}(\boldsymbol{x}) &=&{\frac{1}{nh^d}}\sum_{i=1}^{n}\psi
({\bf Y}_{i})K\Big(
{\frac{\boldsymbol{x}-\boldsymbol{X}_{i}}{h}}\Big);  \\
{\widetilde f}_{{\boldsymbol X},h}({\boldsymbol x})&=&
\E\,f_{\boldsymbol{X},n,h}(\boldsymbol{x})
={\E}\Big\{{\frac{1}{h^d}}K\Big(
{\frac{\boldsymbol{x}-\boldsymbol{X}}{h}}\Big) \Big\}; \\
{\widetilde r}_{\psi,h}({\boldsymbol x})&=&\E\,r_{\psi
,n,h}(\boldsymbol{x}) ={\E}\Big\{{\frac{1}{h^d}}\psi({\bf
Y})K\Big( {\frac{\boldsymbol{x}-\boldsymbol{X}}{h}}\Big) \Big\}.
\end{eqnarray*}
Now, choosing $c_\Psi({\bx})=1/f_{\boldsymbol{X}}({\bx})$ and
$d_\Psi({\bx})=-m_{\Psi }({\bx})/f_{\boldsymbol{X}}({\bx})$ in the
definition $(\ref{proc_emp_multiva_index_fonctions})$ of
$W_{n,h}({\bx},\Psi )$, it is easy to show that, for all $h>0$ the
following relation holds
\begin{equation}
W_{n, h}({\bx},\Psi )={nh^d}\Big\{ \frac{r_{\Psi,n,h} ({\bx})}
{f_{\boldsymbol{X}}({\bx})}-\frac{{\widetilde
r}_{\Psi,h}({\boldsymbol
x})}{f_{\boldsymbol{X}}({\bx})}-\frac{r_{\Psi
}({\bx})\{f_{\boldsymbol{X},{n},h}({\bx})-{\widetilde
f}_{{\boldsymbol X},h}({\boldsymbol x})\}}
{f^2_{\boldsymbol{X}}({\bx})}\Big\} . \label{s1}
\end{equation}
We can now state Lemma \ref{lem_approx_estim_proc_emp_index},
according to which Theorem \ref{theo_reg} is a direct consequence
of Theorem \ref{theo_W}.
\begin{lem}\label{lem_approx_estim_proc_emp_index}
Under the assumptions of theorem \ref{theo_W}, we have, almost
surely,
\begin{eqnarray}
&&{\cal E}_n:= \suph \sup_{\Psi\in\FF_q} \sup_{{\bx}\in
I}\frac{1}{\sqrt{ nh^d\log (1/h^d)}} \Big|W_{n,h}({\bx},\Psi)\nonumber\\[0.2cm]
&&\hskip60pt-nh^d \Big( m_{\Psi;n}({\bx},h)-{r_{\Psi;n}
({\bx},h)\over f_{\boldsymbol{X};n} ({\bx},h)}\Big) \Big| =o(1).
\label{s2}
\end{eqnarray}
\end{lem}
\noindent \demo\ The proof is identical to that of Lemma 10 in
\cite{EM00} (see also Lemma 4.7 in \cite{DM04}) and is omitted for
the sake of conciseness. \findem

\subsection{Proofs of corollaries \ref{cor_fdr}, \ref{cor_density} and
\ref{cor_hazard}}\label{section_proof_cor}

Corollary \ref{cor_fdr} being a direct consequence of Theorem
\ref{theo_carbo} with $\FF=\{\1_{[0,t]}, t\leq\tau_0<T_H\}$, details
of its proof are omitted.

\subsubsection{Proof of Corollary \ref{cor_density}}
Using the following notation
\begin{equation*}
f_{h,\ell}(t;{\bx}):=
\frac{\E\Big\{\1_{\{Y\in[t-\frac{\ell}{2},~t+\frac{\ell}{2}]\}}K\Big(
{{\frac{{\bx}-\boldsymbol{X}}{h}}}\Big) \Big\}}{\ell\E\Big\{K\Big(
\frac{{\bx}-\boldsymbol{X}}{h}\Big) \Big\}},
\end{equation*}
we will first show that
\begin{eqnarray}
\lim_{\ninf} \suph \sup_{\ell\in [\ell'_n,\ell''_n]}
\sup_{{\bx}\in I} \sup_{t\leq\tau_0}\Big\{\fcond(t;{\bx}) -
f_{h,\ell}(t;{\bx})\Big\} =0.\label{stochat}
\end{eqnarray}
Since $\ell''_n\rightarrow 0$, and $\ell'_n\geq \frac{
\sqrt{2\log(1/h'^d_n)}}{nh'^d_n} $ we have for any $s>0$

\begin{eqnarray}&&\lim_{\ninf}
\suph \sup_{\ell\in [\ell'_n,\ell''_n]} \sup_{{\bx}\in I}
\sup_{t\leq\tau_0}\Big\{\fcond(t;{\bx}) - f_{h,\ell}(t;\bx)\Big\}
 \nonumber
\\[0.08in]
&&\leq \lim_{\ninf} \suph  \sup_{{\bx}\in I} \frac{\sqrt{nh^d} \pm
\Theta_{n}({\bx})\sup_{t\leq\tau_0}\sup_{0<\ell\leq
s}\ell\Big\{\fcond(t;{\bx}) -
f_{h,\ell}(t;{\bx})\Big\}}{\sqrt{2\log
(1/h^d)}}.\nonumber
\end{eqnarray}
Set $s\in(0,T_H-\tau_0)$. By applying Theorem \ref{theo_carbo} with
$\FF_s=\{\1_{[t-\ell/2,t+\ell/2]},\ t<\tau_0,\  0<\ell\leq s\}$, we
have with probability one
\begin{eqnarray}
&&\lim_{\ninf} \suph  \sup_{{\bx}\in I} \frac{\sqrt{nh^d} \pm
\Theta_{n}({\bx})\sup_{t\leq\tau_0}\sup_{0<\ell\leq
s}\ell\Big\{\fcond(t;{\bx}) -
f_{h,\ell}(t;{\bx})\Big\}}{\sqrt{2\log (1/h^d)}}
 \nonumber
\\[0.08in]
&&= \Bigg\{ \int_{\R^d}K^{2}(\boldsymbol{u})d\boldsymbol{u}\
\sup_{{\bx}\in I} \frac{\Theta^{2}({\bx})\sup_{\psi\in\FF_s}
\sigma^2_\psi({\bx})}{f_{\boldsymbol{X}}({\bx})} \Bigg\}^{1/2}.
\label{ber1}
\end{eqnarray}
But, by $({\bf A})$, $(F.2)$ and $(D)$,  we have
\begin{eqnarray*}
\lim_{s\rightarrow0}\sup_{{\bx}\in I}\sup_{\psi\in\FF_s}
\sigma^2_\psi({\bx})=0,
\end{eqnarray*}
and then the right-hand term of (\ref{ber1}) can be rendered
arbitrary small, which enables to conclude to
(\ref{stochat}).\vskip5pt

\noindent As for the deterministic term, first observe that

\begin{eqnarray*}
&&\lim_{\ninf} \suph \sup_{\ell\in [\ell'_n, \ell''_n]}
\sup_{{\bx}\in I}
\sup_{t\leq\tau_0}\Big\{f_{h,\ell}(t;{\bx})-f(t;\bx)\Big\}
\nonumber
\\[0.08in]
 &=& \lim_{\ninf} \suph \sup_{\ell\in [\ell'_n, \ell''_n]}
\sup_{{\bx}\in
I}\sup_{t\leq\tau_0}\frac{1}{\frac{1}{h^d}\E\Big\{K\Big(
\frac{{\bx}-\boldsymbol{X}}{h}\Big)
\Big\}f_{\bX}(\bx)}\times\nonumber
\\[0.08in]
&& \Bigg\{
\Big(\frac{1}{h^d\ell}\E\Big\{\1_{\{Y\in[t-\frac{\ell}{2},t+\frac{\ell}{2}]\}}K\Big(
{{\frac{{\bx}-\boldsymbol{X}}{h}}}\Big)
\Big\}-f_{\bX,Y}(\bx,t)\Big)f_{\bX}(\bx)\nonumber
\\[0.08in]
&&+\Big(f_{\bX}(\bx)-\frac{1}{h^d}\E\Big\{K\Big(
\frac{{\bx}-\boldsymbol{X}}{h}\Big) \Big\}\Big)f_{\bX,Y}(\bx,t)
\Bigg\} .
\end{eqnarray*}
By $(F.2)$ and $(D)$ it is enough to establish the two following
statements to conclude the proof of Corollary \ref{cor_density},
\begin{eqnarray}\lim_{\ninf} \sup_{\substack{\inth\\\ell\in
[\ell'_n, \ell''_n]}} \sup_{{\bx}\in
I}\sup_{t\leq\tau_0}\Big\{\frac{1}{h^d\ell}\E\Big[\1_{\{Y\in[t-\frac{\ell}{2},~t+\frac{\ell}{2}]\}}K\Big(
{{\frac{{\bx}-\boldsymbol{X}}{h}}}\Big)
\Big]-f_{\bX,Y}(\bx,t)\Big\}=0\label{biais1}
\end{eqnarray}
and
\begin{eqnarray}\lim_{\ninf} \sup_{\inth} \sup_{{\bx}\in
I}\sup_{t\leq\tau_0}\Big\{f_{\bX}(\bx)-\frac{1}{h^d}\E\Big[K\Big(
\frac{{\bx}-\boldsymbol{X}}{h}\Big) \Big]\Big\}=0.\label{biais2}
\end{eqnarray}
But, for any $h>0$, $\ell>0$, $\bx\in I$, $t<\tau_0$ and $n\geq1$ we
have
\begin{eqnarray*}
&&\frac{1}{h^d\ell}\E\Big\{\1_{\{Y\in[t-\frac{\ell}{2},~t+\frac{\ell}{2}]\}}K\Big(
{{\frac{{\bx}-\boldsymbol{X}}{h}}}\Big) \Big\}-f_{\bX,Y}(\bx,t)\\
&&=\int_{\mathbb{R}^d}\frac{1}{\ell}\int_{t-\frac{\ell}{2}}^{t+\frac{\ell}{2}}\big
(f_{\bX,Y}(\boldsymbol{x}-h\boldsymbol{v},y)-f_{\bX,Y}(\boldsymbol{x}-h\boldsymbol{v},t)\big)dy
K\big(\boldsymbol{v}\big)d\boldsymbol{v}\\
&&\phantom{=}+\int_{\mathbb{R}^d}\big(f_{\bX,Y}(\boldsymbol{x}-h\boldsymbol{v},t)-
f_{\bX,Y}(\bx,t)\big)K\big(\boldsymbol{v}\big)d\boldsymbol{v}\end{eqnarray*}
and
\begin{eqnarray*}f_{\bX}(\bx)-\frac{1}{h^d}\E\Big\{K\Big(
\frac{{\bx}-\boldsymbol{X}}{h}\Big)
\Big\}=\int_{\mathbb{R}^d}\big(f_{\bX}(\boldsymbol{x}-h\boldsymbol{v})-f_{\bX}(\bx)\big)K\big(\boldsymbol{v}\big)d\boldsymbol{v}.\end{eqnarray*}
In view of these two last results, and since $h''_n$ and
$\ell''_n$ decrease to 0 as $\ninf$, it is easy to see that
(\ref{biais1}) and (\ref{biais2}) are direct consequences of
$(K.2)$ and $(D)$.

\subsubsection{Proof of corollary \ref{cor_hazard}}
First observe that
\begin{eqnarray*}
&&\lim_{\ninf} \suph \sup_{\ell\in [\ell'_n,\ell''_n]}
\sup_{{\bx}\in I}
\sup_{t\leq\tau_0}\Big\{\widehat{\lambda}_{n,h,\ell}^{\star}(t;{\bx})
-\lambda(t;\bx)\Big\} \\
&&=\lim_{\ninf} \suph \sup_{\ell\in [\ell'_n,\ell''_n]}
\sup_{{\bx}\in I}
\sup_{t\leq\tau_0}\frac{1}{\big(1-F(t;\bx)\big)\big(1-\Fnhhstar\big)}\times\\
&&\Big\{\big(\fcond(t;\bx)-f(t;\bx)\big)\big(1-F(t;\bx)\big)+f(t;\bx)\big(\Fnhhstar(t;{\bx})-F(t;\bx)\big)\Big\}.
\end{eqnarray*}
In view of Corollaries \ref{cor_fdr} and \ref{cor_density}, and
since $F(\tau_0)<1$ and $f(t;\bx)$ is bounded for $t<\tau_0$ and
$\bx\in I$, it is enough to prove that
\begin{eqnarray*}
&&\lim_{\ninf} \suph\sup_{{\bx}\in I}
\sup_{t\leq\tau_0}\Big(F_{h}(t;{\bx})-F(t;\bx)\Big)=0.
\end{eqnarray*}
But, for any $h>0$, $\bx\in I$, $t<\tau_0$ and $n\geq1$ we have
\begin{eqnarray*}
&&F_{h}(t;{\bx})-F(t;\bx)=\frac{1}{\frac{1}{h^d}\E\Big\{K\Big(
\frac{{\bx}-\boldsymbol{X}}{h}\Big)
\Big\}f_{\bX}(\bx)}\Big\{\int_{y<t}\int_{\mathbb{R}^d}K(\boldsymbol{v})\big(f_{\bX,Y}(\bx-h\boldsymbol{v},y)\\
&&-f_{\bX,Y}(\bx;y)\big)d\boldsymbol{v}dy \times
f_{\bX}(\bx)-\int_{y<t}f_{\bX,Y}(\bx,y)dy\big(\frac{1}{h^d}\E\Big\{K\Big(
\frac{{\bx}-\boldsymbol{X}}{h}\Big) \Big\}-f_{\bX}(\bx)\big)\Big\}.
\end{eqnarray*}
By $(F.2)$, $(D)$ and (\ref{biais2}), it is enough to evaluate the
quantity $\int_{y<t}\int_{\mathbb{R}^d}K(\boldsymbol{v})
\big(f_{\bX,Y}(\bx-$ $h\boldsymbol{v},y)-
f_{\bX,Y}(\bx;y)\big)d\boldsymbol{v}dy$. Towards this aim, set
$\epsilon>0$. Since $I$ is compact and $T_F<\infty$, the
assumption $(D)$ implies the existence of a constant $t(\epsilon)$
such that $\forall \bx\in I^{\alpha}$,
$$\int_{y\leq t(\epsilon)}f_{\bX,Y}(\bx,y)dy\leq \frac{\epsilon}{2}.
$$
So, by $(K.2)$ and $(D)$ and by imposing $h\leq
\min(\frac{2\alpha}{\kappa},
\frac{\epsilon}{2\int_{\mathbb{R}^d}\big|\boldsymbol{v}
K(\boldsymbol{v})\big|d\boldsymbol{v}B_d(\tau_0-t(\epsilon))})$,
we have
\begin{eqnarray*}
&&\Big|\int_{y<t}\int_{\mathbb{R}^d}K(\boldsymbol{v})
\big(f_{\bX,Y}(\bx-h\boldsymbol{v},y)-f_{\bX,Y}(\bx;y)\big)d\boldsymbol{v}dy\Big|\\
&&\leq
\int_{\mathbb{R}^d}\big|K(\boldsymbol{v})\big|\bigg\{\Big|\int_{y\leq
t(\epsilon)}\big(f_{\bX,Y}(\bx-h\boldsymbol{v},y)-f_{\bX,Y}(\bx;y)\big)dy\Big|\\
&&\quad + \Big|\int_{y\in[t(\epsilon);t]}\big(f_{\bX,Y}(\bx-
h\boldsymbol{v},y)-f_{\bX,Y}(\bx;y)\Big)dy\Big|\bigg\}d\boldsymbol{v}\\
&&\leq \epsilon
\end{eqnarray*}
which concludes the proof of Corollary \ref{cor_hazard}. \findem

\newcounter{appendice}

\Alph{appendice}
\setcounter{appendice}{1}

\appendix
\section{Appendix}
\setcounter{equation}{0} \setcounter{theo}{0} \setcounter{lem}{0}
\setcounter{cor}{0} \setcounter{rem}{0} \setcounter{prop}{0}
\setcounter{fact}{0}

\noindent The following fact has been stated in \cite{varron2005}
(see also \cite{Mason04} for a similar result).
\begin{fact}\label{factvarron_reg}
Let $\mathcal{F}$ be a pointwise separable class of functions
satisfying
\begin{equation}
\sup_{f\in\mathcal{F}}{\rm Var}\left(f({\bf Y})\right)\leq
\tau^2h, \nonumber
\end{equation}
with $\tau,h>0$. Assume there exist $M,C,\nu>0$ fulfilling, for
all $0<\epsilon<1,$
\begin{eqnarray*}
\mathcal{N}(\epsilon,\mathcal{F})\leq C\epsilon^{-\nu}, \label{cond1propchap3}\\
\sup_{f\in \mathcal{F},{\bf y}\in\mathbb{R}^d}|f({\bf y})|\leq M
\label{cond2propchap3}.
\end{eqnarray*}
Choose $\rho>0$ arbitrarily. Then, there exist a universal
constant $A_2
>0$ and a parameter $D_1(\nu)>0$ depending only upon $\nu$ such that, if $h>0$
satisfies
\begin{eqnarray*}
K_1:=\max\left\{\frac{4M\sqrt{\nu+1}}{\tau},
\frac{M\rho}{\tau^2}\right\}&\leq& \sqrt{\frac{nh}{\log(1/h)}},\label{cond3propchap3}\\
K_2:=\min\left\{\frac{1}{\tau^2M}, \tau^2\right\}&\geq& h
\label{cond4propchap3},
\end{eqnarray*}
then, setting $T_n(g) = \sum_{j=1}^n \left\{g({\bf Y}_j) -
\mathbb{E}(g({\bf Y}))\right\}$ for $g \in \mathcal{F}$, we have
\begin{equation*}
\mathbb{P}\left(\sup_{1\leq m \leq
n}\|T_m(\cdot)\|_{\mathcal{F}}\geq
(\tau+\rho)D_1\sqrt{nh\log(1/h)}\right) \leq 4
\exp\left(-A_2(\frac{\rho}{\tau})^2\log(1/h)\right).
\end{equation*}
\end{fact}

\bibliographystyle{plain}
%\bibliography{BiblioAM,BiblioNZ}

\end{document}